\input amstex
\documentstyle{amsppt}
\magnification=\magstep1
\hsize=5.2in
\vsize=6.8in
\topmatter
 
\centerline  {\bf NON-ORBIT EQUIVALENT ACTIONS OF $\Bbb F_n$}

\vskip 0.2in
\centerline {\rm ADRIAN IOANA}

\address Math Dept 253-37  
CALTECH, Pasadena, CA 91125\endaddress
\email aioana\@caltech.edu \endemail
\topmatter
\vskip 0.4in
ABSTRACT. For any $2\leq n\leq \infty$, we construct a concrete 1-parameter family of non-orbit equivalent actions of the free group $\Bbb F_n$. These actions arise as  diagonal products between a generalized Bernoulli action and the action  $\Bbb F_n\curvearrowright (\Bbb T^2,\lambda^2)$, where $\Bbb F_n$ is seen as a subgroup of SL$_2(\Bbb Z)$.
\endtopmatter
\document

\vskip 0.1in
\head \S{\bf 0. Introduction}\endhead
\vskip 0.1in

Recall that two free ergodic measure preserving  actions $\Gamma\curvearrowright (X,\mu)$ and $\Lambda\curvearrowright (Y,\nu)$ of two countable discrete groups $\Gamma$ and $\Lambda$ on two standard probability spaces $X$ and $Y$ are said to be {\it orbit equivalent} if there exists a probability space isomorphism $\theta:X\rightarrow Y$ such that $\theta(\Gamma x)=\Lambda \theta(x)$, for $\mu$-almost every $x\in X$.

The orbit equivalence theory of measure preserving group actions has been an extremely active area in the past decade.  New, spectacular rigidity results have been generated using tools ranging from
ergodic theory and operator algebras to representation theory (see the surveys [Ga00],[Sh05],[Po07a]).  
Recently, the problem of finding many non-orbit equivalent actions of a fixed non-amenable group $\Gamma$ has attracted a lot of attention.

This question arose in the 1980's when it was shown that any infinite amenable group $\Gamma$ has exactly one free ergodic measure preserving action, up to orbit equivalence--a result proved by Dye in the case $\Gamma$ is abelian ([Dy59]) and by Ornstein-Weiss in general ([OW80], see  [CFW81] for a generalization)--while some non-amenable groups (e.g. SL$_n(\Bbb Z)$, $n\geq 3$) have uncountably many non-orbit equivalent  actions ([BG81],[Z84],

\noindent
[GG88]). 
In recent years, several classes of non-amenable groups have been shown to share this property: property (T) groups ([Hj05]), weakly rigid groups ([Po06b]), non-amenable products of infinite groups ([Po08], see also [MSh06],[Io07a]) and mapping class groups ([Ki07]). 

In the case of the free groups, progress was slow for a while, only 4 non-orbit equivalent actions of $\Bbb F_n$--all concrete--being known in 2002 ([CW80],[Po06],[Hj05]), before Gaboriau-Popa eventually proved the existence of uncountably many such actions ([GP05]). 
The key idea of their approach was to use the fact that the action of SL$_2(\Bbb Z)$ (as well as its restriction to any free subgroup $\Bbb F_n$) on the 2-torus $\Bbb T^2$ is rigid, in the sense of Popa ([Po06]). 

However, since Gaboriau-Popa's proof also uses a separability argument, it only provides an existence result, leaving open the problem of finding specific actions of $\Bbb F_n$, which are not orbit equivalent. This problem has been emphasized in [Po06a, Section 6], where two further examples were produced, raising the number of concrete non-orbit equivalent actions of $\Bbb F_n$ to 6. 

The main result of this paper is the following:

\vskip 0.1in
\proclaim {Theorem} Let $2\leq n\leq\infty$. Fix an embedding $\Bbb F_n\subset$ SL$_2(\Bbb Z)$ and a surjective homomorphism $\pi:\Bbb F_n\rightarrow\Bbb Z$. Denote by $\sigma$ the restriction of the natural action SL$_2(\Bbb Z)\curvearrowright (\Bbb T^2,\lambda^2)$ to $\Bbb F_n$, where $\lambda^2$ is the Haar measure on $\Bbb T^2$.
For every $t\in (0,1)$, define the probability space $(X_t,\mu_t)=(\{0,1\},r_t)^{\Bbb Z}$, where $r_t(\{0\})=t,r_t(\{1\})=1-t$, and let $\beta_t$ be the Bernoulli action of $\Bbb Z$ on $(X_t,\mu_t)$.
\vskip 0.02in
Let $\alpha_t$ denote the diagonal product action of $\Bbb F_n$ on $(X_t\times \Bbb T^2,\mu_t\times\lambda^2)$ given by $$\alpha_t(\gamma)=\beta_t(\pi(\gamma))\times \sigma(\gamma),\forall\gamma\in\Bbb F_n.$$ 
\vskip 0.02in
 Then $\{\alpha_t\}_{t\in (0,\frac{1}{2}]}$ is a 1-parameter family of free ergodic non-orbit equivalent actions of $\Bbb F_n$.
\endproclaim

To put our main result in a better perspective, note that most non-amenable groups for which concrete uncountable families of non-orbit equivalent actions have been constructed admit in fact  many  actions  
which are {\it orbit equivalent superrigid}, i.e. such that their orbit equivalence class remembers the group and the action. Indeed, this is the case for weakly rigid groups ([Po06b],[Po07]), non-amenable products of infinite groups ([Po08]) and mapping class groups ([Ki07]).
For the free groups, such an extreme rigidity phenomenon never occurs. On the contrary, any free ergodic action of $\Bbb F_n$ is orbit equivalent to actions of uncountably many non-isomorphic groups (see 2.27 in [MSh06]).

The proof of the Theorem has two main parts which we now briefly outline. Assume therefore that $\theta=(\theta_1,\theta_2):X_s\times\Bbb T^2\rightarrow X_t\times\Bbb T^2$ is an orbit equivalence between $\alpha_s$ and $\alpha_t$, for some $s<t\in (0,\frac{1}{2}]$. First we prove that $\theta_i$ ''locally'' (i.e. on a set $A_i\subset X_s\times\Bbb T^2$ of positive measure) depends only 
on the $i$-th coordinate, for $i\in\{1,2\}$. 
This is achieved by playing against each other contrasting properties of the actions $\beta_s$ and $\sigma$. Thus, for $i=1$ we use that $\beta_t$ is an action of an amenable group, while $\sigma$ is strongly ergodic (see Lemma 2.3 and Proposition 2.4) and for $i=2$, we use that $\beta_s$ is a Bernoulli action, whereas $\sigma$ is rigid (see Proposition 3.3). 

For the second part, assume for simplicity that $\theta_i$ depends only on the $i$-th coordinate (i.e. $A_i$ has full measure), for $i\in\{1,2\}$. Letting  $w:\Bbb F_n\times (X_s\times\Bbb T^2)\rightarrow\Bbb F_n$ be the cocycle associated with $\theta$, it follows   that  $\chi=\pi\circ w$ depends only on the $\Bbb F_n$-coordinate. Thus, $\chi$ is a homomorphism $\Bbb F_n\rightarrow\Bbb Z$ which satisfies $\theta_1(\gamma x)=\chi(\gamma)\theta_1(x)$, for all $\gamma\in\Bbb F_n$ and almost every $x\in X_s$. This is further used to prove that $\beta_s$ is isomorphic to the restriction ${\beta_t}_{|m\Bbb Z}$, for some $m\geq 1$. 
 In the general case, we first show that after multiplying $\theta$ with a $\Bbb F_n$-valued function one can assume that $\theta_i$ depends only on the $i$-th coordinate and then proceed as above. This argument, applied to a more general situation, is the subject of Section 4. 
Finally, a simple application of entropy gives that $s\geq t$, a contradiction.

Note that our main result holds for any non-amenable group $\Gamma$ which admits both an infinite amenable quotient $\Delta$ that has no non-trivial finite normal subgroup and a free, weakly mixing, strongly ergodic, rigid action $\Gamma\curvearrowright (Y,\nu)$ (Theorem 5.1).  
\vskip 0.05in
 Recently, a combination of results and ideas from [Io07c], [GL07] and [Ep07] has led to a complete quantitative answer to the problem motivating this paper: any non-amenable group $\Gamma$ admits uncountably many free ergodic non-orbit equivalent actions ([Ep07]).
 Note, however, that the question of finding explicit such actions for an arbitrary $\Gamma$ is still open.

\vskip 0.05in
\noindent
{\it Acknowledgments.} I would like to thank Professor Greg Hjorth and Professor Sorin Popa for useful discussions and encouragement.

\vskip 0.2in

\head \S{\bf 1. Preliminaries} \endhead
\vskip 0.1in

In this section we review some of the notions and results that we will later use.
All  groups $\Gamma$ that we consider hereafter are countable discrete, all probability spaces  $(X,\mu)$ are standard (unless specified otherwise) and all actions $\Gamma\curvearrowright (X,\mu)$ are measure preserving. 
\vskip 0.05in
\noindent
{\bf 1.1 Orbit equivalence and cocycles.} 
Assume that $\Gamma\curvearrowright (X,\mu)$ and $\Lambda\curvearrowright (Y,\nu)$ are two free  orbit equivalent actions. Let $\theta:X\rightarrow Y$ be an orbit equivalence, i.e. a probability space isomorphism such that $\theta(\Gamma x)=\Lambda\theta(x)$, for $\mu$-almost every (a.e.) $x\in X$. 
For every $\gamma\in\Gamma$ and $x\in X$,  denote by $w(\gamma,x)$ the unique (by freeness) element of $\Lambda$ such that $\theta(\gamma x)=w(\gamma,x)\theta(x)$. The map $w:\Gamma\times X\rightarrow\Lambda$ is  measurable, satisfies $$w(\gamma_1\gamma_2,x)=w(\gamma_1,\gamma_2x)w(\gamma_2,x),$$ for all $\gamma_1,\gamma_2\in\Gamma$ and a.e. $x\in X$, and is called the {\it Zimmer cocycle} associated with $\theta$. In general, a measurable map $w:\Gamma\times X\rightarrow \Lambda$ verifying the above relation is called  a {\it cocycle}. Two cocycles $w_1,w_2:\Gamma\times X\rightarrow\Lambda$ are said to be {\it cohomologous} (in symbols, $w_1\sim w_2$) if there exists a measurable map $\phi:X\rightarrow\Lambda$ such that $w_1(\gamma,x)=\phi(\gamma x)w_2(\gamma,x)\phi(x)^{-1}$, for all $\gamma\in\Gamma$ and a.e. $x\in X$.

The simplest instance when two actions $\Gamma\curvearrowright (X,\mu)$ and $\Lambda\curvearrowright (Y,\nu)$ are orbit equivalent is when they are {\it conjugate}, i.e. there exist a probability space isomorphism $\theta:X\rightarrow Y$ and a group isomorphism $\delta:\Gamma\rightarrow\Lambda$ such that $\theta(\gamma x)=\delta(\gamma)\theta(x)$, for all $\gamma\in\Gamma$ and a.e. $x\in X$. Moreover, if $\Gamma=\Lambda$ and $\delta$ is the trivial isomorphism, then we say that the $\Gamma$-actions on $X$ and $Y$ are {\it isomorphic.}  
Much of orbit equivalence rigidity theory aims at proving that, for certain classes of actions, orbit equivalence implies conjugacy. In doing so, the analysis of the associated Zimmer cocycle plays an important role. For example, a general principle proved in [Po07, Proposition 5.11] asserts that if the Zimmer cocycle associated with an orbit equivalence between two weakly mixing actions $\Gamma\curvearrowright (X,\mu)$ and $\Lambda\curvearrowright (Y,\nu)$ is cohomologous to a group homomorphism $\delta:\Gamma\rightarrow\Lambda$, 
then the actions must be (virtually) conjugate.

It is thus very useful to have a criterion for a cocycle to be cohomologous to a group homomorphism. The following theorem, due to S. Popa  (see [Po07, Theorem 3.1]), provides such a criterion. Before stating it, recall that an action $\Gamma\curvearrowright (X,\mu)$ is called {\it weakly mixing} if for every finite collection of measurable sets $A_1,A_2,..,A_n\subset X$ and every $\varepsilon>0$, we  can find $\gamma\in \Gamma$ such that $|\mu(A_i\cap \gamma A_j)-\mu(A_i)\mu(A_j)|\leq\varepsilon,$ for all $i,j\in\{1,..,n\}$. Also, the action $\Gamma\curvearrowright (X,\mu)$ is called {\it mixing} if for every measurable sets $A_1,A_2\subset X$ we have that $\lim_{\gamma\rightarrow\infty}|\mu(A_1\cap \gamma A_2)-\mu(A_1)\mu(A_2)|=0$. 

\proclaim{1.1.1 Theorem [Po07]} Let $\Gamma\curvearrowright (X,\mu)$ be a weakly mixing action and  let $\Gamma\curvearrowright (Y,\nu)$ be another action. Let $\Lambda$ be a countable group and let $w:\Gamma\times (X\times Y)\rightarrow \Lambda$ be a cocycle for the diagonal product action of $\Gamma$ on $X\times Y$. 
Denote by $w^l,w^r:\Gamma\times (X\times X\times Y)\rightarrow\Lambda$ the cocycles for the diagonal product action $\Gamma\curvearrowright X\times X\times Y$ given by $w^l(\gamma,x_1,x_2,y)=w(\gamma,x_1,y)$ and $w^r(\gamma, x_1,x_2,y)=w(\gamma,x_2,y),$ for all $\gamma\in\Gamma,x_1,x_2\in X$ and $y\in Y$.
\vskip 0.03in
If $w^l\sim w^r$, then $w$ is cohomologous to a  cocycle which is independent on the $X$-variable.
\endproclaim

\vskip 0.1in
\noindent
{\bf 1.2 The group measure space construction.}
Let $\Gamma\curvearrowright (X,\mu)$ be a measure preserving action of a countable group $\Gamma$ on a standard probability space $(X,\mu)$.  Let $\Cal H=L^2(X,\mu)\overline{\otimes}\ell^2\Gamma$. For every $\gamma\in\Gamma$ and $f\in L^{\infty}(X,\mu)$, define  the operators $u_{\gamma},L_{f}\in\Bbb B(\Cal H)$ by $$u_{\gamma}(g\otimes \delta_{\gamma'})=\gamma(g)\otimes \delta_{\gamma\gamma'},$$ $$L_f(g\otimes\delta_{\gamma'})=fg\otimes\delta_{\gamma'},\forall\gamma'\in\Gamma,\forall g\in L^2(X,\mu),$$ where, as usual, $\gamma(g)=g\circ\gamma^{-1}$.
Since $u_{\gamma}u_{\gamma'}=u_{\gamma\gamma'},u_{\gamma}L_fu_{\gamma}^*=L_{\gamma(f)}$, for all $\gamma,\gamma'\in\Gamma$ and $f\in L^{\infty}(X,\mu)$, the linear span of $\{L_fu_{\gamma}|f\in L^{\infty}(X,\mu),\gamma\in\Gamma\}$ is a $*$-subalgebra of $\Bbb B(\Cal H)$. The strong operator closure of this algebra, denoted $L^{\infty}(X,\mu)\rtimes\Gamma$, is called the {\it group measure space} von Neumann algebra associated with the action $\Gamma\curvearrowright (X,\mu)$ ([MvN36]). The vector state   $\tau(y)=\langle y (1\otimes\delta_{e}),1\otimes\delta_{e}\rangle$ gives a normal faithful  trace on $L^{\infty}(X,\mu)\rtimes\Gamma$, which is therefore a finite von Neumann algebra.
Furthermore, if the action $\Gamma\curvearrowright (X,\mu)$ is free and ergodic, then $L^{\infty}(X,\mu)\rtimes\Gamma$ is a II$_1$ factor and $L^{\infty}(X,\mu)$ is a Cartan subalgebra, i.e. maximal abelian and regular.

Following [FM77], 
  two free ergodic measure preserving actions  $\Gamma\curvearrowright (X,\mu)$ and $\Lambda\curvearrowright (Y,\nu)$   are   orbit equivalent if and only if  the corresponding Cartan subalgebra inclusions are isomorphic, i.e. $$(L^{\infty}(X,\mu)\subset L^{\infty}(X,\mu)\rtimes \Gamma)\simeq (L^{\infty}(X,\mu)\subset L^{\infty}(X,\mu)\rtimes \Lambda).$$ Moreover, if $\theta:X\rightarrow Y$ is an  orbit equivalence between the actions, then the induced  isomorphism of abelian von Neumann algebras
 $\theta^*:L^{\infty}(Y,\nu)\ni f\rightarrow f\circ\theta\in L^{\infty}(X,\mu)$   extends to an isomorphism $\theta^*:L^{\infty}(Y,\nu)\rtimes\Lambda\rightarrow L^{\infty}(X,\mu)\rtimes\Gamma$. 
We next note that a more general statement of this type is true. Recall first that a measurable map  $q:X\rightarrow Y$ between two probability spaces $(X,\mu)$ and $(Y,\nu)$ is called a {\it quotient map} if it is measure preserving and onto.
In this case, the map $q^{*}:L^{\infty}(Y,\nu)\ni f\rightarrow f\circ q\in L^{\infty}(X,\mu)$ is an embedding of abelian von Neumann algebras. 
\vskip 0.05in

\vskip 0.05in
\proclaim {1.2.1 Lemma [Po07]} Let $\Gamma\curvearrowright (X,\mu),\Lambda\curvearrowright (Y,\nu)$ be two free  actions. Assume that $q:X\rightarrow Y$ is a quotient map such that $q(\Gamma x)=\Lambda q(x)$, a.e. $x\in X$. Also, suppose that  $q$ is 1-1 on the $\Gamma$-orbits, i.e. $q_{|\Gamma x}$ is 1-1, a.e. $x\in X$.
\vskip 0.03in
Then the embedding ${q}^{*}:L^{\infty}(Y,\nu)\hookrightarrow L^{\infty}(X,\mu)$ extends to an embedding $q^{*}:L^{\infty}(Y,\nu)\rtimes\Lambda\hookrightarrow L^{\infty}(X,\mu)\rtimes\Gamma$ of von Neumann algebras.
\endproclaim
This lemma is a particular case of Proposition 1.4.3. in [Po07]. Indeed, if we denote by $\Cal R$ and $\Cal S$ the equivalence relations induced by the actions of $\Gamma$ on $X$ and $\Lambda$ on $Y$, respectively, then $q$ is a local OE of $\Cal R$, $\Cal S$, in the sense of Definition 1.4.2. in [Po07]. By 1.4.3. in [Po07], $q^*$ extends to an embedding $L(\Cal S)\hookrightarrow L(\Cal R)$, where $L(\Cal R)$ denotes the von Neumann algebra associated with $\Cal R$ ([FM77]). Finally, just note that since the actions are assumed free,  $L(\Cal R)$ and $L(\Cal S)$ are naturally isomorphic to $L^{\infty}(X,\mu)\rtimes\Gamma$ and $L^{\infty}(Y,\nu)\rtimes\Lambda$, respectively ([FM77]).
\vskip 0.1in
\noindent
{\bf 1.3 The intertwining bimodule technique.}
This technique has been introduced by S. Popa (see [Po06a, Theorem 2.1 and Corollary 2.3]) and is a powerful tool for deducing unitary conjugacy of subalgebras of a finite von Neumann algebra. Here we note a particular form of it, when the ambient algebra  is abelian. For completeness, we give a self-contained ergodic-theoretic proof. First, we introduce some new terminology.
 \vskip 0.1in

 Let $(X,\mu),(Y,\nu), (Z,\rho)$ be standard probability spaces together with two quotient maps $q:(X,\mu)\rightarrow (Y,\nu)$ and $p:(X,\mu)\rightarrow (Z,\rho)$. Since $q$ is measure preserving, we can
disintegrate $\mu=\int_{Y}\mu_y d\nu(y)$, where $\mu_y$ is a Borel probability measure on $X$ with $\mu_y(q^{-1}(\{y\}))=1$, $\nu$-a.e. $y\in Y$. Let $X\times_{Y}X=\{(x_1,x_2)\in X\times X|q(x_1)=q(x_2)\}$ be the {\it fibered product space} endowed with the probability measure $\mu\times_{\nu}\mu=\int_{Y}(\mu_y\times\mu_y) d\nu(y)$. 
\vskip 0.05in
\noindent
{\bf 1.3.1 Definition.}
 We say that {\it $p$ locally factors through $q$}   if the set 
$$S=\{(x_1,x_2)\in X\times_{Y}X|p(x_1)=p(x_2)\}$$ satisfies $(\mu\times_{\nu}\mu)(S)>0$. 
Equivalently, this means that the set $A$ of $y\in Y$ such that $(\mu_y\times\mu_y)(\{(x_1,x_2)\in X\times X|p(x_1)=p(x_2)\})>0$ has $\nu(A)>0$.
\vskip 0.1in
\noindent
{\bf 1.3.2 Remarks.} $(1)$. To justify our terminology, note that {\it $p$ factors through $q$}, i.e. there exists a quotient map $r:(Y,\nu)\rightarrow (Z,\rho)$ such that $p=r\circ q$, if and only if  $S=X\times_{Y} X$, a.e.

$(2)$. Assume that $(X,\mu)=(Y,\nu)\times (W,\eta)$, for some probability space $(W,\eta)$, and that $q$ is the projection of the $Y$-coordinate.
Then  $p$ locally factors through $q$ if and only if the set $\{(y,w_1,w_2)\in Y\times W\times W|p(y,w_1)=p(y,w_2)\}$ has positive measure.
\vskip 0.05in

\proclaim {1.3.3 Lemma [Po06a]}  Let $(X,\mu),(Y,\nu),(Z,\rho)$ be standard probability spaces  together with two quotient maps $q:(X,\mu)\rightarrow (Y,\nu)$ and $p:(X,\mu)\rightarrow (Z,\rho)$. View $L^{\infty}(Y,
\nu)$ and $L^{\infty}(Z,\rho)$ as von Neumann subalgebras of $L^{\infty}(X,\mu)$, via  $q^*$ and $p^*$, respectively, and let $E:L^{\infty}(X,\mu)\rightarrow L^{\infty}(Y,\nu)$ denote the conditional expectation onto $L^{\infty}(Y,\nu)$. Assume that  there exists $a_1,a_2,..,a_n\in L^{\infty}(X,\mu)$ and $C>0$ such that $$\sum_{i=1}^n||E(fa_i)||_2^2\geq C$$  for all $f\in L^{\infty}(Z,\rho)$ with $|f|=1$ a.e. Then $p$ locally factors through $q$. 
\endproclaim
{\it Proof.} We start by denoting $\tilde X=X\times_{Y}X$ and $\tilde\mu=\mu\times_{\nu}\mu$. Also, for every $f\in L^{\infty}(X,\mu)$,  we define $\tilde f\in L^{\infty}(\tilde X,\tilde\mu)$ by $\tilde f(x_1,x_2)=f(x_1)\overline{f(x_2)}$. Then $$\int_{\tilde X}\tilde f(x_1,x_2)d\tilde\mu(x_1,x_2)=\tag 1.3.a$$ $$\int_{Y}(\int_{q^{-1}(\{y\})\times q^{-1}(\{y\})} f(x_1)\overline{f(x_2)}d\mu_y(x_1)d\mu_y(x_2))d\nu(y)=$$ $$\int_{Y}|\int_{q^{-1}(\{y\})}f(x) d\mu_{y}(x)|^2 d\nu(y)=\int_{Y}|E(f)(y)|^2 d\nu(y)=||E(f)||_2^2.$$ 
Now, let $a\in L^{\infty}(\tilde X,\tilde\mu)$ be given by $a(x_1,x_2)=\sum_{i=1}^n\tilde{a}_i(x_1,x_2)$. Using (1.3.a), the inequality in the hypothesis rewrites as $$\int_{\tilde X} \tilde{f}(x_1,x_2)a(x_1,x_2)d\tilde\mu(x_1,x_2)\geq C\tag 1.3.b$$ for all $f\in L^{\infty}(Z,\rho)$ with $|f|=1$ a.e. 

Next, we denote by $K$ the closed  convex hull of the set $\{\tilde f|f\in L^{\infty}(Z,\rho), |f|=1$  a.e.$\}$ inside the Hilbert space $L^2(\tilde X,\tilde \mu)$.  Let $g\in K$ be the unique element of minimal $||.||_2$. Since $K$ is invariant under the $||.||_2$-preserving transformations $K\ni h\rightarrow \tilde fh\in K$,  we deduce that $g=\tilde f g$,  for every $f\in L^{\infty}(Z,\rho)$ with $|f|=1$ a.e. Thus, if $T=\{(x_1,x_2)\in \tilde X|g(x_1,x_2)\not=0\}$, then for all $f\in L^{\infty}(Z,\rho)$ with $|f|=1$ a.e., we have that $f(p(x_1))=f(p(x_2))$, a.e. $(x_1,x_2)\in T$.

Since $Z$ is a standard probability space we can measurably identify it with the torus $\Bbb T$ endowed with its Haar measure. Thus, by applying the above to the identity function $f(z)=z$, we get that $p(x_1)=p(x_2)$, a.e. $(x_1,x_2)\in T$.
Finally, notice that  (1.3.b) implies that $\int_{\tilde X}g a d\tilde\mu\geq C>0$, hence $g\not=0$ and $\tilde\mu(T)>0$. Altogether, we derive that $\tilde\mu(\{(x_1,x_2)\in\tilde X|p(x_1)=p(x_2)\})>0$, or, in other words, $p$ locally factors through $q$.  \hfill$\square$

\vskip 0.2in

\head \S{\bf 2. Relative strong ergodicity.}\endhead
 \vskip 0.1in 

We begin by recalling that a measure preserving  action $\Gamma\curvearrowright (X,\mu)$ of a countable group $\Gamma$ on a standard probability space $(X,\mu)$ is called {\it strongly ergodic} if for every sequence $\{A_n\}_n\subset X$ of measurable sets satisfying $\lim_{n\rightarrow\infty}\mu(A_n\Delta\gamma A_n)=0$, for all $\gamma\in\Gamma$, we can find  sets $B_n\in\{\emptyset,X\}$ such that $\lim_{n\rightarrow\infty}\mu(A_n\Delta B_n)=0$ ([CW80]).
Examples of  strongly ergodic actions include the actions $\Gamma\curvearrowright (\Bbb T^2,\lambda^2)$, where $\Gamma$ is a non-amenable subgroup of $SL_2(\Bbb Z)$ and $\lambda^2$ is the Haar measure on the 2-torus $\Bbb T^2$ (see [Po06b, Corollary 1.6.5]) and the Bernoulli actions $\Gamma\curvearrowright (X,\mu)^{\Gamma}$ of non-amenable groups $\Gamma$ ([Sc81]).

The notion of strong ergodicity has a useful formulation in terms of von Neumann algebras. Let $\omega$ be a free ultrafilter on $\Bbb N$. The ultraproduct algebra $L^{\infty}(X,\mu)^{\omega}$ is defined as $\ell^{\infty}(\Bbb N,L^{\infty}(X,\mu))/\Cal I_{\omega}$, where $\Cal I_{\omega}$ is the ideal of $f=(f_n)\in\ell^{\infty}(\Bbb N,L^{\infty}(X,\mu))$ for which $\tau_{\omega}(|f|^2)=0$, with the trace $\tau_{\omega}$ being given by $\tau_{\omega}(f)=\lim_{n\rightarrow\omega}\int_{X}f_n d\mu$. Notice that a measure preserving action $\Gamma\curvearrowright (X,\mu)$ induces an integral preserving action of $\Gamma$ on $L^{\infty}(X,\mu)$ which in turn lifts to a $\tau_{\omega}$-preserving action of $\Gamma$ on $L^{\infty}(X,\mu)^{\omega}$. In this context, the action $\Gamma\curvearrowright (X,\mu)$ is strong ergodic if and only if $$[L^{\infty}(X,\mu)^{\omega}]^{\Gamma}:=\{f\in L^{\infty}(X,\mu)^{\omega}|\gamma f=f,\forall\gamma\in\Gamma\}=\Bbb C1.$$

 Moreover, if the action we start with is assumed ergodic, then $[L^{\infty}(X,\mu)^{\omega}]^{\Gamma}$ is either equal to $\Bbb C1$ or infinite dimensional. Here we are making use of the following well-known fact: if $\Gamma\curvearrowright (X,\mu)$ is an ergodic, but not strongly ergodic action, then 
for every $c\in (0,1)$ there exists an asymptotically invariant sequence $\{A_n\}_n\subset X$ such that $\mu(A_n)=c$, for all $n$ (for an idea, see the proof of [JSc87, Lemma 2.3]).

Next, we introduce a relative notion of strong ergodicity.
\vskip 0.05in

\noindent
{\bf 2.1 Definition.} Let $\Gamma\curvearrowright (X,\mu)$ be a measure preserving action together with a quotient action $\Gamma\curvearrowright (Y,\nu)$. Let $q:X\rightarrow Y$ be the associated $\Gamma$-equivariant quotient map and, as usual, view $L^{\infty}(Y,\nu)$ as a von Neumann subalgebra of $L^{\infty}(X,\mu)$.
We say that the action $\Gamma\curvearrowright (X,\mu)$ is {\it strongly ergodic relative to} $\Gamma\curvearrowright (Y,\nu)$ if one of the following equivalent conditions holds true: 
\vskip 0.05in
$(i)$  For every sequence $\{A_n\}_n\subset X$ of measurable sets with $\lim_{n\rightarrow\infty}\mu(A_n\Delta\gamma A_n)=0$, for all $\gamma\in\Gamma$, we can find a sequence of measurable sets $\{B_n\}_n\subset Y$ such that $\lim_{n\rightarrow\infty}\mu(A_n\Delta q^{-1}(B_n))=0$.
\vskip 0.05in
$(ii)$ $[L^{\infty}(X,\mu)^{\omega}]^{\Gamma}=[L^{\infty}(Y,\nu)^{\omega}]^{\Gamma}.$
\vskip 0.05in
$(iii)$ If $f_n\in L^{\infty}(X,\mu), ||f_n||_{\infty}\leq 1$ satisfy $\lim_{n\rightarrow\infty}||f_n-\gamma(f_n)||_2=0$, for all $\gamma\in\Gamma$, then $\lim_{n\rightarrow\infty}||f_n-E(f_n)||_2=0,$ where $E:L^{\infty}(X,\mu)\rightarrow L^{\infty}(Y,\nu)$ denotes the conditional expectation onto $L^{\infty}(Y,\nu)$.
\vskip 0.05in
The proof of the equivalence of conditions $(i)-(iii)$ is standard and we leave it to the reader.
\vskip 0.05in
\noindent
{\bf 2.2 Remarks.} $(1)$ An action $\Gamma\curvearrowright (X,\mu)$ is strongly ergodic if and only if is strongly ergodic relative to the trivial action of $\Gamma$ on a one-point set.

$(2)$ A non-trivial example of relative strong ergodicity arises in the following way. Assume that $\Gamma\curvearrowright (Y,\nu)$ is an action with {\it stable spectral gap}, i.e. such that the unitary representation $\Gamma\curvearrowright (L^2(Y,\nu)\ominus\Bbb C1)\overline{\otimes}(L^2(Y,\nu)\ominus\Bbb C1)$ does not weakly contains the trivial representation (see [Po08, Definition 3.1]). 
Then, for any other measure preserving action $\Gamma\curvearrowright (X,\mu)$, the diagonal product action $\Gamma\curvearrowright (X\times Y,\mu\times\nu)$ is strongly ergodic relative to $\Gamma\curvearrowright (X,\mu)$ (by section 3 in [Po08]).

$(3)$ Note in this respect that if $\Gamma$ is a non-amenable subgroup of SL$_2(\Bbb Z)$, then the action $\Gamma\curvearrowright (\Bbb T^2,\lambda^2)$ has stable spectral gap. Indeed, following the discussion before Lemma 1.6.4 in [Po06b], the representation $\pi$ of $\Gamma$ on $L^2(\Bbb T^2,\lambda^2)\ominus\Bbb C1$ is of the form $\oplus_{i}\ell^2(\Gamma/\Gamma_i)$, where $\{\Gamma_i\}_i$ is a family of amenable subgroups of $\Gamma$. It is easy to see that the product representation $\pi\otimes\pi$ must be of the same form. Thus, by [Po06b, Lemma 1.6.4] we get that $\pi\otimes\pi$ does not weakly contain the trivial representation of $\Gamma$.

\proclaim {2.3 Lemma} 
Let $\Gamma$ be a countable  group and suppose that $\Gamma_0\subset\Gamma$ is a normal subgroup such that the quotient group $\Delta=\Gamma/\Gamma_0$ is infinite amenable. Assume that $\Gamma\curvearrowright (X,\mu)$ is a free strongly ergodic measure preserving  action such that its restriction to $\Gamma_0$ is ergodic. Let $\Delta\curvearrowright (Y,\nu)$ be a free ergodic measure preserving action and let $\Gamma$ act on $Y$ via the homomorphism $\Gamma\rightarrow\Delta$. 

Then the diagonal product action $\Gamma\curvearrowright (X\times Y,\mu\times \nu)$ is strongly ergodic relative to the quotient $\Gamma\curvearrowright (Y,\nu)$.
\endproclaim
{\it Proof.}  We first show that since the action $\Gamma\curvearrowright (X,\mu)$ is strongly ergodic, its restriction to $\Gamma_0$ must also be strongly ergodic. If we assume the contrary, then $\Cal A:=[L^{\infty}(X,\mu)^{\omega}]^{\Gamma_0}\not=\Bbb C1.$ Moreover, as the action $\Gamma_0\curvearrowright (X,\mu)$ is ergodic, $\Cal A$ must be infinite dimensional. 
Since $\Gamma_0$ is a normal subgroup of $\Gamma$, we get that $\Gamma$ acts on $\Cal A$ and that this action passes to an action of $\Delta$. 

Towards a contradiction, we claim that $ ({\Cal A}^{\omega})^{\Delta}\not=\Bbb C1.$ Remark that we can  assume that $\Cal A^{\Delta}=\Bbb C1$, since otherwise the claim follows trivially. 
Next,
let $\Cal B$ be an infinite dimensional, $\Delta$-invariant, separable von Neumann subalgebra of $\Cal A$. 
 Since $\Cal B$ is abelian we can identify it with the $L^{\infty}$-algebra of a probability space $(Z,\rho)$ in such a way that the action of $\Delta$ on $\Cal B$ is induced by a measure preserving action $\Delta\curvearrowright (Z,\rho)$. The fact that $\Cal B$ is infinite dimensional  implies that $(Z,\rho)$ is not completely atomic. Moreover, since the action $\Delta\curvearrowright (Z,\rho)$ is ergodic (as $\Cal B^{\Delta}=\Bbb C1$), it follows that $(Z,\rho)$ is  completely non-atomic.

Since $\Delta$ is amenable, we can apply Ornstein-Weiss' theorem ([OW80]) to derive that the action $\Delta\curvearrowright (Z,\rho)$ is orbit equivalent to a free ergodic action of $\Bbb Z$, and thus is not strongly ergodic. Therefore, we get that $(\Cal B^{\omega})^{\Delta}\not=\Bbb C1$ and furthermore that $ ({\Cal A}^{\omega})^{\Delta}\not=\Bbb C1$.
This however implies that $[L^{\infty}(X,\mu)^{\omega}]^{\Gamma}\not=\Bbb C1$, or, equivalently, that the action $\Gamma\curvearrowright (X,\mu)$ is not strongly ergodic, a contradiction.
\vskip 0.05in
Now, to prove the conclusion of the lemma, let $\{A_n\}_n\subset X\times Y$ be a sequence of measurable sets such that $\lim_{n\rightarrow\infty}(\mu\times\nu)(A_n\Delta\gamma A_n)=0$, for all $\gamma\in\Gamma$. For every $n$ and $y\in Y$, denote $A_n^y=\{x\in X|(x,y)\in A_n\}$.
Since $\Gamma_0$ acts trivially on $Y$ we have that
$$(\mu\times\nu)(A_n\Delta\gamma A_n)=\int_{Y}\mu(A_n^y\Delta\gamma A_n^y)d\nu(y),\forall n\geq 0,\gamma\in\Gamma_0.$$ 

As by our assumption  $\lim_{n\rightarrow\infty}(\mu\times\nu)(A_n\Delta\gamma A_n)=0$, for all $\gamma\in\Gamma_0,$
then, after eventually passing to a subsequence of $\{A_n\}$, we may assume that $$\lim_{n\rightarrow\infty}\mu(A_n^y\Delta\gamma A_n^y)=0,\forall\gamma\in\Gamma_0,$$ a.e. $y\in Y$. The strong ergodicity of the action $\Gamma_0\curvearrowright (X,\mu)$  implies that $$\lim_{n\rightarrow\infty}\mu(A_n^y)(1-\mu(A_n^y))=0,$$ a.e. $y\in Y$. Thus, if we denote $B_n=\{y\in Y|\lim_{n\rightarrow\infty}\mu(A_n^y)=1\}$, then $\lim_{n\rightarrow\infty}(\mu\times\nu)(A_n\Delta (X\times B_n))=0$, which proves the lemma.
 \hfill $\square$
\vskip 0.1in

Jones and Schmidt showed that an ergodic action $\Gamma\curvearrowright (X,\mu)$ is  strongly ergodic if and only if given any free ergodic $\Sigma\curvearrowright (Z,\rho)$ of an infinite amenable group $\Sigma$ there does not exist a quotient map $p
:X\rightarrow Z$ such that $p(\Gamma x)=\Sigma p(x)$, a.e. $x\in X$ ([JSc87]). Next, we generalize the only if part of their result to a relative strong ergodicity situation. We will later use this generalization to analyze the orbit equivalences between certain diagonal product actions (see the proof of Theorem 4.1).

\proclaim {2.4 Proposition}  Let $\Gamma\curvearrowright (X,\mu)$ be a measure preserving action which is strongly ergodic relative to a quotient action $\Gamma\curvearrowright (Y,\nu)$. Let  $q:X\rightarrow Y$ be the associated $\Gamma$-equivariant quotient map. Assume that $\Sigma$ is an infinite amenable group and let $\Sigma\curvearrowright (Z,\rho)$ be a free ergodic measure preserving action. 

If $p:X\rightarrow Z$ is a quotient map such that $p(\Gamma x)=\Sigma p(x)$, a.e. $x\in X$, then $p$ locally  factors through $q$. \endproclaim

{\it Proof.} 
 By Ornstein-Weiss' theorem ([OW80]), the action 
 $\Sigma\curvearrowright (Z,\rho)$ is orbit equivalent to the diagonal action  $(\oplus_{m\geq 0}\Bbb Z_2)\curvearrowright \prod_{m\geq 0}(\{0,1\},r)_m$, where $r$ is the probability measure on $\{0,1\}$ with both weights equal to $\frac{1}{2}$ and $\Bbb Z_2$ acts on $\{0,1\}$ in the only non-trivial way.
We can therefore assume that $(Z,\rho)=\prod_{m\geq 0}(\{0,1\},r)_m$ and that $\Sigma z=(\oplus_{m\geq 0}\Bbb Z_2) z,$ a.e. $z\in Z$. For $n\geq 0$, we define $(Z_n,\rho_n)=\prod_{m\geq n}(\{0,1\},r)_m$ and view $L^{\infty}(Z_n,\rho_n)$ as a von Neumann subalgebra of $L^{\infty}(Z,\rho)$, via the projection $\pi_n:Z\rightarrow Z_n$.
\vskip 0.05in
\noindent
{\it Claim 1.} For every $n$, let $f_n\in L^{\infty}(Z_n,\rho_n)$ with $||f_n||_{\infty}\leq 1$. Then  $$\lim_{n\rightarrow\infty}||f_n\circ p-\gamma(f_n\circ p)||_2=0,\forall\gamma\in\Gamma.$$
\vskip 0.05in
{\it Proof of Claim 1.} For each $\gamma\in\Gamma$ and every $n\geq 1$, define the measurable set $C_{\gamma,n}=\{x\in X|p(\gamma^{-1}x)\in (\oplus_{m<n} \Bbb Z_2) p(x)\}$. Since $p(\Gamma x)=\Sigma p(x)=(\oplus_{m\geq 0}\Bbb Z_2)p(x)$, a.e. $x\in X$, we deduce that $\lim_{n\rightarrow\infty}\mu(C_{\gamma,n})=1$, for all $\gamma\in\Gamma$.  Now, for all $n$ and for every $f\in L^{\infty}(Z_n,\rho_n)$, we have that $$f(z)=f(gz), \forall z\in Z,\forall g\in (\oplus_{m<n}\Bbb Z_2).$$
Thus $$(\gamma (f_n\circ p))(x)=f_n(p(\gamma^{-1}x))=f_n(p(x))=(f_n\circ p)(x),$$ for all $x\in C_{\gamma,n}$ and since  $||f_n||_{\infty}\leq 1$, we get that $||f_n\circ p-\gamma (f_n\circ p)||_2\leq 2\sqrt{1-\mu
(C_{\gamma,n})},$ which proves the claim.
\vskip 0.05in
\noindent
{\it Claim 2.} There exists $n$ such that for all $f_n\in L^{\infty}(Z_n,\rho_n)$ with $||f_n||_{\infty}\leq 1$ we have that $$||f_n\circ p-E(f_n\circ p)||_2\leq \frac{1}{2},$$ where $E:L^{\infty}(X,\mu)\rightarrow L^{\infty}(Y,\nu)$ denotes the conditional expectation onto $L^{\infty}(Y,\nu)$.
\vskip 0.05in
{\it Proof of Claim 2.} Assuming the claim false, for every $n$ we can find $f_n\in L^{\infty}(Z_n,\rho_n)$ with $||f_n||_{\infty}\leq 1$ and $||f_n\circ p-E(f_n\circ p)||_2>\frac{1}{2}.$ If we define $f=(f_n\circ p)_n\in L^{\infty}(X,\mu)^{\omega}$, then by Claim 1 we get that $f\in [L^{\infty}(X,\mu)^{\omega}]^{\Gamma}$. By the relative strong ergodicity assumption we deduce that $f\in [L^{\infty}(Y,\nu)]^{\Gamma}$ and thus $\lim_{n\rightarrow\omega}||f_n\circ p-E(f_n\circ p)||_2=0,$ which contradicts the above.
\vskip 0.05in
\noindent
{\it Claim 3.} $p$ locally factors through $q$.
\vskip 0.05in
{\it Proof of Claim 3.} The previous claim implies  that $$||E(f_n\circ p)||_2\geq \frac{1}{2},$$ for all $f_n\in L^{\infty}(Z_n,\rho_n),||f_n||_{\infty}\leq 1$ with $|f_n|=1$ a.e. 
By applying Lemma 1.3.3 we deduce  that $\pi_n\circ p:X\rightarrow Z_n$ locally factors through $q$. Thus, if we disintegrate $\mu=\int_{Y}\mu_y d\nu(y)$, then there is a set $A\subset Y$ with $\nu(A)>0$ such  that $(\mu_y\times\mu_y)(\{(x_1,x_2)\in X\times X|(\pi_n\circ p)(x_1)=(\pi_n\circ p)(x_2)\})>0$, for all $y\in A$. 

Using Fubini's theorem, for every $y\in A$, we can find $z_y\in Z_n$ such that $\mu_y(\{x\in X|(\pi_n\circ p)(x)=z_y\})>0$. Further, if we let $\pi^n:Z\rightarrow\prod_{0\leq m<n}\{0,1\}_m$ be the projection onto the first $n$-coordinates, then for every $y\in A$, we can find $t_y\in\prod_{0\leq m<n}\{0,1\}_m$ such that  $$\mu_y(\{x\in X|p(x)=(t_y,z_y)\})=\mu_y(\{x\in X|(\pi_n\circ p)(x)=z_y,(\pi^n\circ p)(x)=t_y\})>0.$$

 Finally, the last inequality implies that $(\mu_y\times\mu_y)(\{(x_1,x_2)\in X\times X|p(x_1)=p(x_2)\})>0$, for all $y\in A$, and thus $p$ locally factors through $q$.
\hfill$\square$

\vskip 0.2in

\head \S{\bf 3. Rigidity and Bernoulli actions.}\endhead
\vskip 0.1in

We first review S. Popa's notion of rigidity for actions.

\vskip 0.05in
 \noindent
{\bf 3.1 Rigid actions.}
 Let $N$ be a separable, finite von Neumann algebra together with a faithful, normal trace $\tau$ and let $B\subset N$ be a von Neumann subalgebra. 
On $N$ we consider the $2$-norm given by $||x||_2=\tau(x^*x)^{1/2}$, for all $x\in N$. The inclusion $B\subset N$ is called {\it rigid} (cf. [Po06, section 4])  if for any  sequence $\phi_n:N\rightarrow N$ of unital, tracial, completely positive maps such that $\phi_n\rightarrow 1_N$ (i.e. $\lim_{n\rightarrow \infty}||\phi_n(x)-x||_2=0$, for all $x\in N$) we must have that $$\lim_{n\rightarrow \infty}\sup_{x\in B,||x||\leq 1}||\phi_n(x)-x||_2=0.$$ 

A free ergodic measure preserving action $\Gamma\curvearrowright (X,\mu)$ is called {\it rigid} if its associated Cartan subalgebra inclusion $L^{\infty}(X,\mu)\subset L^{\infty}(X,\mu)\rtimes \Gamma$ is rigid. For example, if $\Gamma$ is a non-amenable subgroup of $\Bbb F_2$, then the action of $\Gamma\curvearrowright (\Bbb T^2,\lambda^2)$ is rigid, where $\lambda^2$ denotes the Haar measure of the $2$-torus $\Bbb T^2$. This follows by combining the fact (proven in [Bu91]) that the pair $(\Gamma\ltimes\Bbb Z^2,\Bbb Z^2)$ has the relative property (T) of Kazhdan-Margulis with Proposition 5.1 in [Po06].

\vskip 0.1in

Next, we briefly recall the definition of generalized Bernoulli actions.
\vskip 0.05in
\noindent
{\bf 3.2 Bernoulli actions.}
 Let $\Gamma$ be a countable group, $I$ be a countable set on which $\Gamma$ acts (e.g. $I=\Gamma/\Gamma_0$, for some subgroup $\Gamma_0$ of $\Gamma$)  and $(X_0,\mu_0)$ be a probability space. The measure preserving action $\Gamma\curvearrowright (X_0,\mu_0)^I$ given by $\gamma((x_i)_i)=((x_{{\gamma}^{-1}\cdot i})_i)$, for all $x=(x_i)_i\in (X_0,\mu_0)^I$ and each $\gamma\in\Gamma$, is called a {\it generalized Bernoulli action}.
In the case $I=\Gamma$, with $\Gamma$ acting on itself by left multiplication,  we will call such an action a {\it Bernoulli action}.
\vskip 0.1in

In recent years, S. Popa proved remarkable rigidity results concerning Bernoulli actions (see the survey [Po07a]). The general philosophy behind these results is that any rigidity phenomenon that is exhibited by the group measure space, $M$, associated with a Bernoulli action $\Gamma\curvearrowright (X,\mu)=(X_0,\mu_0)^{\Gamma}$, has to come, in some sense, from the group $\Gamma$. 
Following this principle, it is proven in [Io07b] (see Corollary 3.7 for $B=L^{\infty}(X_0,\mu_0)$) that if $A\subset M$ is a  rigid inclusion of von Neumann algebras, then $A$ can be essentially conjugated, via a unitary, into the group von Neumann algebra $L\Gamma$. In particular, one derives that $A$ cannot be a diffuse subalgebra of $L^{\infty}(X,\mu)$.

A direct consequence of this result is the following fact: given any rigid action $\Lambda\curvearrowright (Z,\rho)$ of an infinite group $\Lambda$, there is no quotient map $p:X\rightarrow Z$ such that  $p(\Gamma x)=\Lambda p(x)$ and $p_{|\Gamma x}$ is 1-1, a.e. $x\in X$. Indeed, if there is such a $p$, then Lemma 1.2.1 would imply that the inclusion $A=p^*(L^{\infty}(Z,\rho))\subset M$  is rigid. This is, however, a contradiction since $A\subset L^{\infty}(X,\mu)$.

More generally, we have:
\proclaim {3.3 Proposition} Let $\Gamma\curvearrowright I$ be an action of a countable group $\Gamma$ on a countable set $I$ and let $(X_0,\mu_0)$ be a (possibly atomic) probability space. Denote $(X,\mu)=(X_0,\mu_0)^I$ and let $\Gamma\curvearrowright (X,\mu)$ be the generalized Bernoulli action  induced by the action of $\Gamma$ on $I$. Let $\Gamma\curvearrowright (Y,\nu)$ be another measure preserving action and on $X\times Y$ consider the diagonal product action of $\Gamma$. Denote by $q:X\times Y\rightarrow Y$ the projection onto the $Y$-coordinate.

Let $\Lambda\curvearrowright (Z,\rho)$ be a free ergodic rigid measure preserving action of an infinite countable group $\Lambda$.  Suppose that $p:X\times Y\rightarrow Z$ is a quotient map such that $p(\Gamma(x,y))=\Lambda p(x,y)$ and $p_{|\Gamma(x,y)}$ is 1-1, a.e. $(x,y)\in X\times Y$.
\vskip 0.03in
Then $p$ locally factors through $q$.
\endproclaim

{\it Proof.} We begin by encoding the hypothesis in the language of von Neumann algebras. Let $p^*:L^{\infty}(Z,\rho)\rightarrow L^{\infty}(X\times Y,\mu\times\nu)$ be the embedding given by $p^*(f)=f\circ p$. Using the hypothesis, Lemma 1.2.1 implies that $p^*$ extends to an embedding $p^*:L^{\infty}(Z,\rho)\rtimes\Lambda\rightarrow 
  M:=L^{\infty}(X\times Y,\mu\times\nu)\rtimes\Gamma$. Moreover, since the inclusion $L^{\infty}(Z,\rho)\subset L^{\infty}(Z,\rho)\rtimes\Lambda$ is assumed rigid, we get that the inclusion $A:=p^*(L^{\infty}(Z,\rho))\subset M$ is also rigid (see 4.6. in [Po06]).
 
\vskip 0.1in
Before continuing, we need to review some  von Neumann algebra notions. 
\vskip 0.1in
\noindent
{\bf 3.4 Terminology}. Hereafter, we will work with finite von Neumann algebras $N$ endowed with a fixed faithful normal trace $\tau$.  
We denote by $L^2(N)$  the completion  of $N$ with respect to the norm $||.||_2$. The scalar product on $L^2(N)$ therefore verifies  $\langle x,y\rangle=\tau(y^*x)$, for all $x,y\in N$. 
Note that $N$ is also endowed with the operator norm $||.||$ and that $||xy||_2\leq \min\{||x|| ||y||_2,||x||_2 ||y||\}$, for all $x,y\in N$.
 Let $P$ be a von Neumann subalgebra of $N$, which we will always assume to be endowed with the trace $\tau_{|P}$. The orthogonal projection from $L^2(N)$ onto $L^2(P)$ takes $N$ onto $P$ and its restriction to $N$ is precisely the unique $\tau$-preserving conditional expectation, $E_P$, from $N$ onto $P$. Recall that $E_P$ is $P$-bimodular, i.e. $E_P(p_1xp_2)=p_1E_P(x)p_2$, for all $p_1,p_2\in P$ and $x\in N$, and that $E_P$ is a contraction in both of the above norms.
We also denote by $\Cal U(N)$ the group of unitaries of $N$, i.e. elements $u\in N$ such that $u^*u=1$. 
Every unitary element $u$ induces an automorphism Ad$(u)$ of $N$ through the formula Ad$(u)(x)=uxu^*$. 

Next, we briefly recall three fundamental constructions involving von Neumann algebras.
 Given an (always assumed $\tau$-preserving) action $\alpha:\Gamma\rightarrow$ Aut$(N)$ of a countable group $\Gamma$ on $N$, the associated {\it crossed product von Neumann algebra} $N\rtimes_{\alpha} \Gamma$ is defined in the same way as the group measure space algebra $L^{\infty}(X,\mu)\rtimes\Gamma$, where one replaces $L^{\infty}(X,\mu)$ by $N$ throughout the construction in 1.2. Note that every element $x\in N\rtimes_{\alpha}\Gamma$ can be uniquely written as $x=\sum_{\gamma\in\Gamma}x_{\gamma}u_{\gamma}$, where $x_{\gamma}\in N$, for all $\gamma\in\Gamma$. The trace $\tau$ on $N$ extends to a trace on  $\tilde\tau$ on $N\rtimes_{\alpha}\Gamma$ through the formula $\tilde\tau(x)=\tau(x_e)$.

If $(N_1,\tau_1)$ and $(N_2,\tau_2)$ are two finite von Neumann algebras, then $N_1\overline{\otimes}N_2$ denotes their {\it tensor product von Neumann algebra} endowed with the trace $\tau=\tau_1\otimes\tau_2$. For example, if $(N_i,\tau_i)=(L^{\infty}(X_i,\mu_i),\int d\mu_i)$, for some probability spaces $(X_i,\mu_i)$ ($i=1,2$), then $(N_1\overline{\otimes} N_2,\tau_1\otimes\tau_2)$ is naturally isomorphic to $(L^{\infty}(X_1\times X_2,\mu_1\times\mu_2),\int d(\mu_2\times\mu_2))$. Now, if $\theta_i$ is automorphism of $N_i$, then $\theta_1\otimes\theta_2$ denotes the automorphism of $N_1\overline{\otimes} N_2$ given by $(\theta_1\otimes\theta_2)(x_1\otimes x_2)=\theta_1(x_1)\otimes\theta_2(x_2)$, for all $x_i\in N_i$. Thus, the diagonal product action $\alpha=\alpha_1\times\alpha_2:\Gamma\rightarrow$ Aut$(N_1\overline{\otimes}N_2)$ of two actions $\alpha_i:\Gamma\rightarrow$ Aut$(N_i)$ is defined by $\alpha(\gamma)=\alpha_1(\gamma)\otimes\alpha_2(\gamma)$.
Given a von Neumann algebra $N$ and a countable set $I$, we denote by $N^I$ the tensor product von Neumann algebra $\overline{\otimes}_{i\in I}(N)_i$. If $J$ is a subset of $I$, then we view $N^J$ as a subalgebra of $N^I$, via the isomorphism $N^J\cong (\overline{\otimes}_{i\in J}(N)_i)\overline{\otimes}(\overline{\otimes}_{i\in I\setminus J}(\Bbb C1)_i)$.

Finally, the {\it free product} of two finite von Neumann algebras $(N_1,\tau_1)$ and $(N_2,\tau_2)$ is the unique finite von Neumann algebra $(N,\tau)$ (in symbols, $N=N_1*N_2$) such that: (1) $N$ contains both $N_1$ and $N_2$, (2) $\tau_{|N_i}=\tau_i$, for $i\in\{1,2\}$, (3) $N_1$ and $N_2$ generate $N$ as a von Neumann algebra and (4) $\tau(x_{i_1}x_{i_2}..x_{i_n})=0$, for all $n\geq 1$, $i_1\not=i_2\not=..\not= i_n\in\{1,2\}$ and all $x_{i_k}\in N_{i_k}$ with $\tau_{i_k}(x_{i_k})=0$, for all $k$.  
\vskip 0.1in
Going back to the proof of Proposition 3.3, we define $B=L^{\infty}(X_0,\mu_0)$ and $C=L^{\infty}(Y,\nu)$. 
Then $ L^{\infty}(X,\mu)\cong B^I$ and the action $\Gamma\curvearrowright (X,\mu)$ induces a Bernoulli action $\alpha:\Gamma\rightarrow$ Aut$(B^I)$. Similarly, let $\beta:\Gamma\rightarrow$ Aut$(C)$ be the action induced by $\Gamma\curvearrowright (Y,\nu)$. With these notations, $M$ can be viewed as the crossed product von Neumann algebra $(B^I\overline{\otimes} C)\rtimes_{(\alpha\times\beta)}\Gamma$. To summarize, at this point, we know that $A$ is a von Neumann subalgebra of $B^I\overline{\otimes} C$ such that the inclusion $A\subset M$ is rigid.
\vskip 0.05in

Following Popa's deformation/rigidity strategy we will use the deformability properties of Bernoulli actions against the rigidity of the inclusion $A\subset M$ in
 order to determine the position of $A$ inside $M$. 

\vskip 0.1in
\noindent
{\it Step 1.} First, we recall a construction from [Io07b, Proposition 2.3]. More precisely,  we augment $M$ to a von Neumann algebra $\tilde M$ and define a $1$-parameter group of automorphisms $\{\Theta_t\}_{t\in\Bbb R}$ of $\tilde M$ such that $\Theta_t\rightarrow 1_{\tilde M}$, as $t\rightarrow 0$. Towards this, we define  the free product von Neumann algebra $\tilde B=B*L^{\infty}(\Bbb T,\lambda)$, where $\lambda$ is Haar measure on the torus $\Bbb T$.
 Let $\tilde\alpha:\Gamma\rightarrow \text{Aut}({\tilde B}^I)$ be the Bernoulli action given by the action of $\Gamma$ on $I$. It is clear that  $B^I\subset {\tilde B}^I$ and that $\tilde\alpha$ extends $\alpha$, hence we have the inclusion $$M\subset \tilde M:=({\tilde B}^I\overline{\otimes} C)\rtimes_{(\tilde\alpha\times\beta)}\Gamma.$$
 We denote by $\tau$ the natural trace on $\tilde M$, i.e. the trace obtained from the integration traces on $B$ and $C$ by applying the corresponding constructions from 3.4.

Now, let $u\in L^{\infty}(\Bbb T,\lambda)$ be the Haar unitary given by $u(z)=z$, for all $z\in\Bbb T$, and let $h\in L^{\infty}(\Bbb T,\lambda)$ be a real-valued function such that $u=e^{ih}$. For every $t\in\Bbb R$, we define the unitary element  $u_t=e^{ith}\in L^{\infty}(\Bbb T,\lambda)$ (thus $u_t\in\tilde B$) and consider the automorphism $$\theta_t=(\otimes_{i\in I}\text{Ad}(u_t)_i)\otimes 1_C\in \text{Aut}({\tilde B}^I\overline{\otimes}C)\tag 3.a$$

Notice that  $\theta_t$ commutes with the action $\tilde\alpha\times\beta$, for all $t$.  Thus, $\theta_t$ extends to an automorphism $\Theta_t\in\text{Aut}(\tilde M)$ given by  $$\Theta_t(x)=\sum_{\gamma\in\Gamma}\theta_t(x_{\gamma})u_{\gamma},$$ for all $x=\sum_{\gamma\in\Gamma}x_{\gamma}u_{\gamma}\in \tilde M$, where $x_{\gamma}\in {\tilde B}^I\overline{\otimes} C,$ for all $\gamma\in\Gamma$. Since  $\lim_{t\rightarrow 0}||u_t-1||_2=0$, we get that $\theta_t\rightarrow 1_{{\tilde B}^I\overline{\otimes} C}$,  thus $\Theta_t\rightarrow 1_{\tilde M}$, as $t$ goes to 0.
\vskip 0.1in
\noindent
{\it Step 2.} Secondly, we use the rigidity of the inclusion $A\subset M$ to deduce that for some $t>0$ we can find a non-zero $v\in {\tilde B}^I\overline{\otimes} C$ such that $\theta_t(u)v=vu$, for all $u\in\Cal U(A)$. Indeed, 
since the inclusion $A\subset M$ is rigid, we get that the inclusion $A\subset\tilde M$ is rigid (by  [Po06, 4.6]). Thus, since $\Theta_t\rightarrow 1_{\tilde M}$, we can find $t>0$ such that $$||\theta_t(x)-x||_2=||\Theta_t(x)-x||_2\leq \frac{1}{2},\forall x\in A,||x||\leq 1.$$ 

 We next employ a standard functional analysis trick (see e.g. the proof of 4.4 in [Po06]). Let $\Cal K$ be the $||.||_2$-closed convex hull of the set $\{\theta_t(u)u^*|u\in\Cal U(A)\}$ and note that if $k\in \Cal K$ and $u\in\Cal U(A)$, then $\theta_t(u)ku^*\in \Cal K$ and $||\theta_t(u)ku^*||_2=||k||_2$. Thus, if $v$ denotes the unique element of minimal $||.||_2$ in $\Cal K$, then $v\in {\tilde B}^I\overline{\otimes} C$, $||v||\leq 1$ and $\theta_t(u)vu^*=v$, for all $u\in\Cal U(A)$. Thus $\theta_t(u)v=vu$, for all $u\in\Cal U(A)$.  Moreover, since $||\theta_t(u)u^*-1||_2=||\theta_t(u)-u||_2\leq \frac{1}{2}$, for all unitaries $u\in A$, we get that $||v-1||_2\leq \frac{1}{2},$ hence $v\not=0$.
\vskip 0.1in

We denote by $\Cal H$ the Hilbert space $L^2({\tilde B}^I\overline{\otimes} C)$. For every subset $F$ of $I$, we let $P_F$  be the orthogonal projection of $\Cal H$ onto $L^2(B^F\overline{\otimes}C)$. Also,  we denote by $Q_F$ and $R_F$ the orthogonal projections of $\Cal H$ onto $L^2({\tilde B}^F\overline{\otimes}C)$ and onto $L^2(\tilde{B}^F\overline{\otimes}{(u_t B u_t^*)}^{I\setminus F}\overline{\otimes}C)$, respectively.

\vskip 0.1in
\noindent
{\it Step 3.} In this context, we next prove that we can find $c>0$ and $F\subset I$ finite such that  $$||P_F(u)||_2\geq c,\forall u\in\Cal U(A)\tag 3.b$$

\vskip 0.05in
To this end, let $t>0$ and $v$ be as given by Step 2 and assume for simplicity that $||v||\leq 1$.  
We start by approximating $v$ with a finitely supported vector $w$. 
Notice that if $F_n$ are increasing finite subsets of $I$ such that $\cup_{n}F_n=I$, then $Q_{F_n}\rightarrow 1_{L^2({\tilde B}^I\overline{\otimes}C)}$, in the strong operator topology. Thus, for large enough $n$,  $w=Q_{F_n}(v)$ satisfies $$c=(1-|\tau(u_t)|^2)||w||_2-2||w-v||_2>0\tag 3.c$$ (here we are using the fact that $|\tau(u_t)|<1$, for all $t>0$). Set $F=F_n$ and note that $w\in {\tilde B}^F\overline{\otimes}C$ and $||w||\leq||v||\leq 1$.

Now, if we fix $u\in\Cal U(A)$, then since $u\in B^I\overline{\otimes}C$ and $w\in {\tilde B}^F\overline{\otimes}C$, we get that $\theta_t(u)v\in L^2({\tilde B}^F\overline{\otimes}(u_tBu_t^*)^{I\setminus F}\overline{\otimes}C)$ and thus that $R_F(\theta_t(u)w)=\theta_t(u)w$. 
By combining this fact with the equality $\theta_t(u)v=vu$ and using triangle's inequality, we  derive that $$||R_F(wu)-wu||_2=||R_F(wu-\theta_t(u)w)-(wu-\theta_t(u)w)||_2\leq\tag 3.d$$ $$2||wu-\theta_t(u)w||_2=2||(w-v)u-\theta_t(u)(w-v)||_2\leq 2||w-v||_2.$$
Further, (3.d) and (3.c) together imply that $$||R_F(wu)||_2\geq ||wu||_2-2||w-u||_2=||w||_2-2||w-u||_2\tag 3.e$$ $$c+|\tau(u_t)|^2||w||_2,\forall u\in \Cal U(A).$$ 
\noindent
Next, we estimate from above the expression $||R_F(wu)||_2$.
\proclaim {3.5 Lemma} For all $\omega\in L^2({\tilde B}^F\overline{\otimes} B^{I\setminus F}\overline{\otimes}C)$, we have that 

$$||R_F(\omega)||_2^2\leq |\tau(u_t)|^4||\omega||_2^2 + (1-|\tau(u_t)|^4)||Q_F(\omega)||_2^2.$$
\endproclaim
\noindent
We postpone the proof of this lemma until the end of this section. 
For the moment, we assume the lemma true and explain how it finishes the proof of Step 3. Indeed, if $u\in\Cal U(A)$, then it is clear that $wu\in {\tilde B}^F\overline{\otimes} B^{I\setminus F}\overline{\otimes}C$. 
Since $Q_F$ is ${\tilde B}^F\overline{\otimes}C$-bimodular and $w\in{\tilde B}^F\overline{\otimes}C$, we get that $Q_F(wu)=wQ_F(u)$. Also, since $u\in B^I\overline{\otimes}C$, it is clear that $Q_F(u)=P_F(u)$. Thus, $||Q_F(wu)||_2=||wP_F(u)||_2\leq ||w|| ||P_F(u)||_2\leq ||P_F(u)||_2$.

Altogether, we can apply Lemma 3.5 to $\omega=wu$ to deduce that $$||R_F(wu)||_2^2\leq |\tau(u_t)|^4||wu||_2^2 + (1-|\tau(u_t)|^4)||Q_F(wu)||_2^2\leq$$ $$|\tau(u_t)|^4||w||_2^2 + (1-|\tau(u_t)|^4)||P_F(u)||_2^2,\forall u\in\Cal U(A).$$ It is now clear that this inequality combined with (3.e) proves the claim of Step 3.
\vskip 0.1in
\noindent
{\it Step 4.} We can now finish the proof of Proposition 3.3. 
First, we strengthen the conclusion of Step 3 by showing that there exists a finite dimensional subalgebra $D\subset B^I$ such that $$||E_{D\otimes C}(u)||_2\geq c/2,\forall u\in\Cal U(A)\tag 3.f$$ where $E_{D\otimes C}$ denotes the conditional expectation onto $D\otimes C$. To see this, let $B_n \subset B$ be an increasing sequence of finite dimensional $*$-algebras such that the $*$-algebra $\cup_{n\geq 1}B_n$ is dense in $B$, in the strong operator topology. Using the rigidity of the inclusion $A\subset M$, the same argument as in the proof of 5.3.1. in [Po06] shows that there exists $n$ such that $||E_{B_n^I\overline{\otimes}C}(u)-u||_2\leq c/2$, for all $u\in\Cal U(A)$. Indeed, this follows by  just noticing that if $E_n$ denotes the conditional expectation onto $M_n=(B_n^I\overline{\otimes}C)\rtimes_{(\alpha\times\beta)}\Gamma$, then  $E_n\rightarrow 1_{M}$ and $E_n(u)=E_{B_n^I\overline{\otimes}C}(u)$, for all $u\in\Cal U(A)$. 

Since $P_F$ is a contraction, we further get that $||P_F(E_{B_n^I\overline{\otimes}C}(u))-P_F(u)||_2\leq c/2$, for every unitary $u\in A$. On the other hand, by Step 3,  $||P_F(u)||_2\geq c$. Therefore combining the last two inequality yields $$||(P_F\circ E_{B_n^I\overline{\otimes}C})(u)||_2\geq c/2,\forall u\in\Cal U(A).$$ It is clear that if $D=B_n^F$, then $D$ is finite dimensional and $P_F\circ E_{B_n^I\overline{\otimes}C}=E_{D\otimes C}$. Thus, (3.f) is proven.

Next, since $D$ is a finite dimensional abelian algebra, we can find projections $p_1,..,p_n$ such that $D=\Bbb Cp_1\oplus...\oplus\Bbb Cp_n$. A simple calculation then shows that $E_{D\otimes C}(x)=\sum_{i=1}^n(p_i\otimes E_C(xp_i))/\tau(p_i)$, for all $x\in B^I\overline{\otimes}C.$ Thus, $$||E_{D\otimes C}(x)||_2^2=\sum_{i=1}^n||E_C(xp_i)||_2^2/\tau(p_i),\forall x\in B^I\overline{\otimes}C.$$
From this identity and (3.f) it follows that $\sum_{i=1}^n||E_C(up_i)||_2^2/\tau(p_i)\geq c/2,$ for all $u\in\Cal U(A)$. Finally, after noticing that for every $f\in L^{\infty}(Z,\rho)$ with $|f|=1$, a.e., $u=f\circ p$ is a unitary in $A$, we can apply Lemma 1.3.3 to deduce the conclusion.\hfill$\square$

\vskip 0.1in

In the proof of Lemma 3.5 we will need the following result:
\proclaim {3.6 Lemma}  If $P_1$  denotes the orthogonal projection of $L^2(\tilde B)$ onto $L^2(u_tBu_t^*)$, then $P_1(1)=1$ and $P_1(\zeta)=|\tau(u_t)|^2(u_t\zeta u_t^*),$ for all $\zeta\in L^2B\ominus\Bbb C1.$ 
\endproclaim
{\it Proof.}
Using the $||.||_2$-density of $B$ in $L^2B$ we can assume that $\zeta\in B$ with $\tau(\zeta)=0$. Fix $b\in B$ and set $b_0=b-\tau(b)$. Also, let $z_t=u_t-\tau(u_t)$. Then we have that $$\langle P_1(\zeta),u_tbu_t^*\rangle=\langle\zeta,u_tbu_t^*\rangle=\langle \zeta,u_tb_0u_t^*\rangle=\tau(u_t^*b_0^*u_t\zeta)=$$ $$\tau(z_t^*b_0^*z_t\zeta)+\tau(u_t)\tau(z_t^*b_0^*\zeta)+\overline{\tau(u_t)}\tau(b_0^*z_t\zeta)+|\tau(u_t)|^2\tau(b_0^*\zeta).$$

Since $\zeta,b_0\in B$ and $z_t\in L^{\infty}(\Bbb T,\lambda)$ all have zero traces, by using the fact that $B$ and $L^{\infty}(\Bbb T,\lambda)$ are in a free position, we deduce that the first three terms in the last sum are equal to 0. Thus, we have that $$ \langle P_1(\zeta),u_tbu_t^*\rangle=|\tau(u_t)|^2\tau(b_0^*\zeta)=|\tau(u_t)|^2\langle\zeta,b_0\rangle=$$ $$|\tau(u_t)|^2\langle \zeta,b\rangle=\langle |\tau(u_t)|^2 (u_t\zeta u_t^*),u_tbu_t^*\rangle,\forall b\in B,$$
 which, again by the density of the inclusion $B\subset L^2B$, implies the conclusion of the lemma.\hfill$\square$
\vskip 0.1in
{\it Proof of Lemma 3.5.} If  $\Cal K=L^2({\tilde B}^F\overline{\otimes
}C)$, then we identify $\Cal H=L^2({\tilde B}^I\overline{\otimes} C)$ with the tensor product Hilbert space $(\overline{\otimes}_{i\in I\setminus F}L^2\tilde B)\overline{\otimes} \Cal K,$ in the natural way.  Under this identification, we have that $R_F=(\otimes_{i\in I\setminus F}P_1)\otimes 1_{\Cal K}$, while $Q_F=
 (\otimes_{i\in I\setminus F}P_0)\otimes 1_{\Cal K},$ where $P_0$ denotes the orthogonal projection of $L^2\tilde B$ onto $\Bbb C1$ (in other words, $P_0(x)=\tau(x)$, for all $x\in\tilde B$). 

Next, using this identification, we construct a specific orthonormal basis $\Omega$ for the Hilbert space $L^2({\tilde B}^F\overline{\otimes} B^{I\setminus F}\overline{\otimes}C)$ which is now identified with $(\overline{\otimes}_{i\in I\setminus F}L^2B)\overline{\otimes} \Cal K.$ To this end, let  $\{\xi_j\}_{j\geq 0}$, and $\{\eta_k\}_{k\geq 0}$ be orthonormal bases for  the Hilbert spaces $L^2B$ and $\Cal K$, respectively, and assume that $\xi_0=\eta_0=1$. Let $S$ be the set of pairs $(j,k)\in \Bbb N^{I\setminus F}\times\Bbb N$  such that $S_j:=\{i\in I\setminus F|j_i\not=0\}$ is finite. For every $(j,k)\in S$, we define $$\omega_{(j,k)}=(\otimes_{i\in I\setminus F}\xi_{j_i})\otimes \eta_k.$$ Then clearly  $\Omega:=\{\omega_{(j,k)}|(j,k)\in S\}$ is an orthonormal basis for $(\overline{\otimes}_{i\in I\setminus F}L^2B)\overline{\otimes} \Cal K.$ 

By Pythagoras' theorem in order to prove the lemma it is  sufficient to check that
(1) the conclusion of the lemma holds for every $\omega=\omega_{(j,k)}\in\Omega$ and that (2)  $R_F(\omega)\perp R_F(\omega')$ and $Q_F(\omega)\perp Q_F(\omega')$,  for every $\omega\not=\omega'\in\Omega$. Indeed, using the  tensor product decomposition of $R_F$ and Lemma 3.6, we have that $$R_F(\omega_{(j,k)})=(\otimes_{i\in I\setminus F}P_1(\xi_{j_i}))\otimes \eta_k=|\tau(u_t)|^{2|S_j|}(\otimes_{i\in S_j}u_t\xi_{j_i}u_t^*)\otimes \eta_k.$$ Similarly, $Q_F(\omega_{(j,k)})$ is equal to $\eta_k$, if $S_j=\emptyset$, and to 0, otherwise. Using these formulas, (2) is immediate. To check (1), we differentiate two cases. If $S_j=\emptyset$, then $R_F(\omega)=Q_F(\omega)=\omega$ (all being equal to $\eta_k$), while if $S_j\not=\emptyset$, then $||R_F(\omega)||_2=|\tau(u_t)|^{2|S_j|}\leq |\tau(u_t)|^2||\omega||_2$. In conclusion, (1) holds true in both cases.\hfill$\square$

\vskip 0.2in
\head \S{\bf 4. Orbit equivalence and product actions}\endhead
\vskip 0.1in

\proclaim{ 4.1 Theorem}
Let $\Gamma,\Lambda$ be two countable  groups and suppose that:
\vskip 0.02in
 $(i)$ $\Delta=\Gamma/\Gamma_0$ and $\Sigma=\Lambda/\Lambda_0$ are infinite amenable quotients  of $\Gamma$ and $\Lambda$.
 
 $(ii)$ $\Delta$ has no non-trivial finite normal subgroup.

$(iii)$  $\Delta\curvearrowright (X_1,\mu_1)=(X_1^0,\mu_1^0)^{\Delta}$ is a Bernoulli action.

$(iv)$ $\Gamma\curvearrowright(X_2,\mu_2)$ is a free, weakly mixing, strongly ergodic, measure preserving action. 

$(v)$ $\Sigma\curvearrowright(Y_1,\nu_1)$ is a free, mixing, measure preserving action.

$(vi)$ $ \Lambda\curvearrowright(Y_2,\nu_2)$ is a free, ergodic, rigid, measure preserving action.
\vskip 0.02in
  Let $\Gamma$ act on $X_1$ and $\Lambda$ act on $Y_1$ through the homomorphisms $\Gamma\rightarrow \Delta$ and $\Lambda\rightarrow \Sigma$. 
\vskip 0.05in
Assume that the diagonal product actions $\Gamma\curvearrowright (X_1\times X_2,\mu_1\times\mu_2)$ and $\Lambda\curvearrowright (Y_1\times Y_2,\nu_1\times\nu_2)$ are orbit equivalent.  Then a quotient $\Delta\curvearrowright (X_0,\mu_0)$ of the action $\Delta\curvearrowright (X_1,\mu_1)$  is conjugate to the restriction  of the action $\Sigma\curvearrowright (Y_1,\nu_1)$  to some subgroup $\Sigma_0\subset\Sigma$.
\endproclaim
{\it Proof.} 
\vskip 0.05in
\noindent
{\it Step 1.} Let $\theta=(\theta_1,\theta_2):X_1\times X_2\rightarrow Y_1\times Y_2$ be a measure space isomorphism such that $\theta(\Gamma x)=\Lambda\theta(x)$, a.e. $x\in X_1\times X_2$. In the first part of the proof we use the results of the previous two sections to show that $\theta_i$ ''locally'' only depends on the $i$-th coordinate. To this end, let $q_i:X_1\times X_2\rightarrow X_i$ denote the projection onto the $i$-th coordinate, for $i\in\{1,2\}$.

Firstly, since the action $\Gamma\curvearrowright (X_2,\mu_2)$ is strongly ergodic and $\Delta$ is amenable, by Lemma 2.3 we get that the action $\Gamma\curvearrowright (X_1\times X_2,\mu_1\times \mu_2)$ is strongly ergodic relative to the quotient $\Gamma\curvearrowright (X_1,\mu_1)$. Since we also have that $\theta_1(\Gamma x)=\Lambda\theta_1(x)=\Sigma\theta_1(x)$, a.e. $x\in X_1\times X_2$, and $\Sigma$ is amenable, we can thus use Proposition 2.4 to deduce that $\theta_1$ locally factors through $q_1$.

Secondly, note that $\theta_2(\Gamma x)=\Lambda\theta_2(x)$ and ${\theta_2}_{|\Gamma x}$ is 1-1, a.e. $x\in X_1\times X_2$. Since $\Gamma\curvearrowright (X_1,\mu_1)$ is a generalized Bernoulli action, while the action $\Lambda\curvearrowright (Y_2,\nu_2)$ is assumed rigid, we are in a position to apply Proposition 3.3 and thereby get that $\theta_2$ locally factors through $q_2$.

Altogether, in view of Remark 1.3.2 (2), we conclude that the sets $\{(x_1,x_2,x_2')\in X_1\times X_2\times X_2|\theta_1(x_1,x_2)=\theta_1(x_1,x_2')\}$ and $\{(x_1,x_1',x_2)\in X_1\times X_1\times X_2|\theta_2(x_1,x_2)=\theta_2(x_1',x_2)\}$ both have positive measures. 
Let us exploit this information further and show that $$\theta_1(x_1,x_2)\in\Lambda\theta_1(x_1,x_2'),\theta_2(x_1,x_2)\in\Lambda\theta_2(x_1',x_2)\tag 4.a$$ a.e. $(x_1,x_1',x_2,x_2')\in X_1\times X_1\times X_2\times X_2.$
Indeed, from the above we get that the set $T=\{(x_1,x_2,x_2')\in X_1\times X_2\times X_2|\theta_1(x_1,x_2)\in\Lambda\theta_1(x_1,x_2')\}$ has positive measure. On the other hand, $T$ is clearly invariant under the diagonal product action of $\Gamma$ on $X_1\times X_2\times X_2$. Since this action is weakly mixing (as all actions in the product are) it follows that $T=X_1\times X_2\times X_2$, a.e. This proves the first assertion in (4.a) and the second one follows similarly.

\vskip 0.1in
\noindent
{\it Step 2.}
Next, we use (4.a) in combination with S. Popa's criterion for untwisting cocycles (Theorem 1.1.1) to show that some perturbation $\rho=\phi$ $\theta$ of $\theta$ with a function $\phi:X\rightarrow\Lambda$ is of the form $\rho=(\rho_1,\rho_2)$, where $\rho_i:X_i\rightarrow Y_i$ is a measurable map, for $i\in\{1,2\}$. Moreover, we prove that $\rho_1$ verifies the formula $\rho_1(\gamma x_1)=\chi(\gamma)\rho_1(x_1)$, for some group homomorphism $\chi:\Gamma\rightarrow\Sigma$.
To start with, let $w:\Gamma\times (X_1\times X_2)\rightarrow\Lambda$ be the Zimmer cocycle associated with $\theta$.

By (4.a), we can find a measurable map $\psi_2:X_1\times X_1\times X_2\rightarrow\Lambda$ such that $\theta_2(x_1,x_2)=\psi_2(x_1,x_1',x_2)\theta_2(x_1',x_2)$, a.e. $(x_1,x_1',x_2)\in X_1\times X_1\times X_2$. Using the fact that $\Lambda$ acts freely on $Y_2$, we deduce that $$\psi_2(\gamma x_1,\gamma x_1',\gamma x_2)w(\gamma,(x_1',x_2))=w(\gamma,(x_1,x_2))\psi_2(x_1,x_1',x_2)\tag 4.b$$ for all $\gamma\in\Gamma$ and a.e. Indeed, it is clear that if we apply both sides of (4.b) to $\theta_2(x_1',x_2)$ and use the definitions of $w$ and $\psi_2$ then we obtain the same result. 

Since the action $\Gamma\curvearrowright (X_1,\mu_1)$ is weakly mixing, it follows from Theorem 1.1.1 that we can find a cocycle $w_2:\Gamma\times X_2\rightarrow\Lambda$ and a measurable map $\phi_2:X_1\times X_2\rightarrow\Lambda$ such that 
$w(\gamma,(x_1,x_2))=\phi_2(\gamma x_1,\gamma x_2)^{-1}w_2(\gamma,x_2)\phi_2(x_1,x_2)$, for all $\gamma\in\Gamma$ and a.e. $(x_1,x_2)\in X_1\times X_2$. Combining this with the definition of $w$ further yields that $$\phi_2(\gamma x_1,\gamma x_2)\theta_i(\gamma x_1,\gamma x_2)=w_2(\gamma,x_2)\phi_2(x_1,x_2)\theta_i(x_1,x_2)\tag 4.c$$ for each $\gamma\in\Gamma$, a.e. $(x_1,x_2)\in X_1\times X_2$ and all $i\in\{1,2\}$.

Now, we denote $\eta=\phi_2$ $\theta:X_1\times X_2\rightarrow Y_1\times Y_2$ and claim that if $\eta=(\eta_1,\eta_2)$, then  $\eta_2(x_1,x_2)=\eta_2(x_1',x_2)$, a.e. $(x_1,x_1',x_2)\in X_1\times X_1\times X_2$. Let $U$ denote the set of such triples $(x_1,x_1',x_2)$. By (4.c) it is clear that $U$ is invariant under the diagonal $\Gamma$-action and since this action is ergodic (being a product of weakly mixing actions) it is therefore sufficient to show that $U$ has positive measure. In turn, this is a consequence of the following three facts: $\eta_2=\phi_2$ $\theta_2$, $\phi_2$ takes countably many values and the set $\{(x_1,x_1',x_2)\in X_1\times X_1\times X_2|\theta_2(x_1,x_2)=\theta_2(x_1',x_2)\}$ has positive measure.
\vskip 0.1in
To summarize, at this point we have that a perturbation $\eta=\phi_2$ $\theta$ of $\theta$ is of the form $\eta=(\eta_1,\eta_2)$, with $\eta_2:X_2\rightarrow Y_2$ (by the claim), and satisfies $\eta(\gamma x_1,\gamma x_2)=w_2(\gamma,x_2)\eta(x_1,x_2)$, for every $\gamma\in\Gamma$ and a.e. $(x_1,x_2)\in X_1\times X_2$ (by (4.c)).

Our plan is now to repeat the above arguments, with $\eta_1$ instead of $\theta_2$. Note first that by (4.a), $\eta_1(x_1,x_2)\in \Lambda\eta_1(x_1,x_2')$, a.e. $(x_1,x_2,x_2')\in X_1\times X_2\times X_2$. Thus, we can find a measurable map $\psi_1:X_1\times X_2\times X_2\rightarrow\Lambda$ such that  $\eta_1(x_1,x_2)=\psi_1(x_1,x_2,x_2')\eta_1(x_1,x_2')$, a.e.
 Let $\pi:\Lambda\rightarrow\Sigma$ be the quotient homomorphism. By taking into account that $\Sigma$ acts freely on $Y_1$, we deduce that $$\pi(\psi_1(\gamma x_1,\gamma x_2,\gamma x_2'))\pi(w_2(\gamma,x_2'))=\pi(w_2(\gamma,x_2))\pi(\psi_1(x_1,x_2,x_2')),$$ for all $\gamma\in\Gamma $ and a.e. $(x_1,x_2,x_2')\in X_1\times X_2\times X_2$. As before, one verifies this formula by applying both sides to $\eta_1(x_1,x_2')$. Since the action $\Gamma\curvearrowright (X_2,\mu_2)$ is weakly mixing, Theorem 1.1.1  gives a measurable map $\phi_1:X_2\rightarrow\Lambda$ and a group homomorphism $\chi:\Gamma\rightarrow\Sigma$ such that $\pi(\phi_1(\gamma x_2)w(\gamma,x_2)\phi_1(x_2)^{-1})=\chi(\gamma),$ for each $\gamma\in\Gamma$ and a.e. $x_2\in X_2$.
Thus, if we define $\rho=(\rho_1,\rho_2):X_1\times X_2\rightarrow Y_1\times Y_2$ by $\rho(x_1,x_2)=\phi_1(x_2)\eta(x_1,x_2)$, then for a.e. $(x_1,x_2)\in X_1\times X_2$ and all $\gamma\in\Gamma$ we have that $$\rho_1(\gamma x_1,\gamma x_2)=\pi(\phi_1(\gamma x_2))\eta_1(\gamma x_1,\gamma x_2)=\pi(\phi_1(\gamma x_2)w_2(\gamma,x_2))\eta_1(x_1,x_2)=$$
$$\chi(\gamma)\pi(\phi_1(x_2))\eta_1(x_1,x_2)=\chi(\gamma)\rho_1(x_1,x_2)$$
The same argument as above now shows that $\rho_1$ only depends on the $X_1$-coordinate. Thus, the above formula rewrites as $$\rho_1(\gamma x_1)=\chi(\gamma)\rho_1(x_1)\tag 4.d$$for each $\gamma\in\Gamma$ and a.e. $x_1\in X_1$.  It is also clear by construction that $\rho_2$ only depends on the $X_2$-coordinate. To complete the proof of this step, just notice that if $\phi:X_1\times X_2\rightarrow\Lambda$ is given by $\phi(x_1,x_2)=\phi_1(x_2)\phi_2(x_1,x_2)$, then $\rho=\phi$ $\theta$.
\vskip 0.1in
\noindent
{\it Step 3.} In this  step we aim to show that $\rho_1$ is a quotient map, $\chi(\Gamma_0)=e$ and that the induced group homomorphism $\chi:\Delta=\Gamma/\Gamma_0\rightarrow\Sigma$ is injective. Prior to proving these assertions, let us indicate how they imply the conclusion of the theorem.
Indeed, assuming them true, let $(X_0,\mu_0)=(Y_1,\nu_1)$ and note that $\Delta$ acts on $X_0$ by $\gamma x_0=\chi(\gamma) x_0$, for each $\gamma\in\Gamma$ and all $x_0\in X_0$. It is now clear that $\rho_1:X_1\rightarrow X_0$ is a $\Delta$-equivariant (by (4.d))  quotient map  and that the identity map $X_0\rightarrow Y_1$ conjugates the actions $\Delta\curvearrowright (X_0,\mu_0)$ and $\Sigma_0=:\chi(\Delta)\curvearrowright (Y_1,\nu_1)$.
\vskip 0.05in
Turning to the proof of the above assertions, we first claim that there exists a measurable set $A_1\subset X_1$ such that $\mu_1(A_1)>0$, $\nu_1(\rho_1(A_1))>0$ and $\rho_1:(A_1,c{\mu_1}_{|A_1})\rightarrow (\rho_1(A_1),{\nu_1}_{|\rho_1(A)})$ is a measure space isomorphism, where $c=\frac{\nu_1(\rho_1(A_1))}{\mu_1(A_1)}$. To see this, for every $\gamma\in\Gamma$, let $X_{\gamma}=\{x\in X_1\times X_2|\phi(x)=\gamma\}$. Then $\{X_{\gamma}\}_{\gamma\in\Gamma}$ is a measurable partition of $X_1\times X_2$, so, in particular, we can find $\gamma$ such that $(\mu_1\times\mu_2)(X_{\gamma})>0$. After discarding a measure zero set from $X_{\gamma}$ we can assume that $\theta_{|X_{\gamma}}$ is 1-1 and since $\rho=\gamma$ $\theta$ on $X_{\gamma}$, we deduce that $\rho_{|X_{\gamma}}$ is 1-1.
This implies that there exist measurable sets $A_i\subset X_i$ such that $\mu_i(A_i)>0$ and ${\rho_i}_{|A_i}$ is 1-1, for $i\in\{1,2\}$. Thus, $\rho_{|A}$ is 1-1, where $A=A_1\times A_2$.
Moreover, note that $\rho_{|X_{\gamma}\cap A}:X_{\gamma}\cap A\rightarrow \rho(X_{\gamma}\cap A)$ is  measure preserving (being equal to $\gamma$ $\theta$), for all $\gamma\in\Gamma$. Thus, we get that $\rho_{|A}:A\rightarrow\rho(A)$ is a measure preserving isomorphism (where on $A$ and $\rho(A)$ we consider the restrictions of the measures $\mu_1\times\mu_2$ and $\nu_1\times\nu_2$, respectively) and since $\rho=(\rho_1,\rho_2)$, the claim follows. 

Next, we prove that $\chi(\Gamma_0)=e$. This is a consequence of the following three facts:  (1) $\chi(\gamma)$ stabilizes $\rho_1(x_1)$, a.e. $x_1\in X_1$ and for all $\gamma\in\Gamma_0$ (by (4.d)), (2) $\nu_1(\rho_1(C))>0$ for every $C\subset X_1$ such that $\mu_1(X_1\setminus C)=0$ (by the above claim) and (3) $\Sigma$ acts freely on $Y_1$. Now, let $\chi$ also denote the induced homomorphism $\Delta=\Gamma/\Gamma_0\rightarrow\Sigma$. Since $\Delta$ has no non-trivial finite normal subgroups, in order to prove that $\chi$ is injective, it is  enough to show that Ker$(\chi)$ is finite. Assume by contradiction that Ker$(\chi)$ is infinite and let $B\subset Y_1$ be a measurable set. By (4.d) we get that $\rho_1^{-1}(B)$ is Ker$(\chi)$-invariant and since the action $\Delta\curvearrowright (X_1,\mu_1)$ is mixing (hence its restriction to Ker$(\chi)$ is ergodic), we would derive that $\mu_1(\rho_1^{-1}(B))\in\{0,1\}$. On the other hand, if we choose $B\subset\rho_1(A_1)$ such that $\nu_1(B)\in (0,\nu_1(\rho_1(A_1)))$, then, by using the claim proved above, it is clear that  $\mu_1(\rho_1^{-1}(B))\in (0,1),$ a contradiction.

Finally, let us show that $\rho_1$ is a quotient map. By (4.d) we have that $\rho_1(X_1)$ is $\Sigma_0$-invariant, while the claim insures that $\nu_1(\rho_1(X_1))>0$. Since $\Sigma_0\cong\Delta$ is infinite and the action $\Sigma\curvearrowright (Y_1,\nu_1)$ is mixing, we get that the restriction $\Sigma_0\curvearrowright (Y_1,\nu_1)$ is ergodic, hence it follows that $\rho_1(X_1)=Y_1$. Thus, it remains to show that $\rho_1$ is measure preserving. In other words, we need to prove that if $\nu$ is the probability measure on $Y_1$ given by $\nu(B)=\mu_1(\rho_1^{-1}(B))$, for every measurable subset $B\subset Y_1$, then $\nu=\nu_1$.
 Notice that by (4.d) we get that $\gamma(\rho_1^{-1}(B))=\rho_1^{-1}(\chi(\gamma)B)$, for each $\gamma\in\Delta$ and every subset $B$ of $Y_1$. This implies that $\nu$ is  $\Sigma_0$-invariant. Since the action $\Sigma_0\curvearrowright (Y_1,\nu_1)$ is ergodic, to finish the proof, it is thus sufficient to show that $\nu$ is absolutely continuous with respect to $\nu_1$. Let $B\subset Y_1$ such that $\nu_1(B)=0$. Then for every $\gamma\in\Delta$ we have that 

$$\mu_1(\rho_1^{-1}(B)\cap \gamma^{-1}A_1)=\mu_1(\gamma\rho_1^{-1}(B)\cap A_1)=\mu_1(\rho_1^{-1}(\chi(\gamma)B)\cap A_1)$$ which by the above claim is further equal to $c^{-1}\nu_1(\chi(\gamma)B\cap \rho_1(A_1))$ and thus to 0 (since $\nu_1(\chi(\gamma)B)=\nu_1(B)=0$). Using the ergodicity of the action $\Delta\curvearrowright (X_1,\mu_1)$ we get that $X_1=\cup_{\gamma\in\Delta}\gamma A_1$, which altogether implies that $\mu_1(\rho_1^{-1}(B))=0$, as needed.
\hfill$\square$

\vskip 0.2in
\head \S{\bf 5. Proof of the main result}\endhead
\vskip 0.1in
In this section, we derive the main result as a consequence of Theorem 4.1.
Before stating and proving a more general version of the main result (Theorem 5.1), let us recall a few well-known facts about entropy (see [OW87] and [Pe83] for a reference). 
Let $\Delta$ be an infinite amenable group. Given  a measure preserving action $\sigma:\Delta\curvearrowright (X,\mu)$ of $\Delta$ on a standard probability space $(X,\mu)$   we denote by $h(\sigma)$ its entropy. If $\sigma_0:\Delta\curvearrowright (X_0,\mu_0)$ is a quotient action of $\sigma$, then $h(\sigma_0)\leq h(\sigma)$. In particular, two isomorphic actions have the same entropy.

For $n\geq 1$ and a $n$-tuple $p=(p_1,p_2,..,p_n)$ of positive numbers with sum equal to $1$,
 let $(X_p,\mu_p)$ be the product probability space $(\{1,2,..,n\},r_p)^{\Delta},$ where $r_p(\{i\})=p_i$, for all $i\in\{1,2,..,n\}$. 
If $\beta_p$ denotes the Bernoulli action of $\Delta$ on $(X_p,\mu_p)$, then $h(\beta_p)=-\sum_{i=1}^{n}p_i\log_2(p_i).$
 On the other hand, if $(Z,\rho)$ is a standard probability space then the entropy of the Bernoulli action $\Delta\curvearrowright (Z,r)^{\Delta}$ is equal to $+\infty$.
Given a subgroup $\Delta_0$ of $\Delta$, the restriction ${\beta_p}_{|\Delta_0}$ is precisely the Bernoulli  action of $\Delta_0$ with base  $(\{1,2,..,n\},r_p)^{\Delta/\Delta_0}$. Thus, using the above remarks, it is easy to see that $h({\beta_p}_{|\Delta_0})=|\Delta/\Delta_0|h(\beta_p).$

\vskip 0.1in
\proclaim {5.1 Theorem } Let $\Gamma$ be a countable group and assume that 
\vskip 0.02in
$(i)$ $\Delta$ is an infinite amenable quotient of $\Gamma$ together with a surjective homomorphism $\pi:\Gamma\rightarrow\Delta$,
 \vskip 0.02in
$(ii)$ $\Delta$ has no non-trivial finite normal subgroup and
\vskip 0.02in
$(iii)$ $\sigma:\Gamma\curvearrowright (Y,\nu)$ is a free, weakly mixing, strongly ergodic, rigid measure preserving action of $\Gamma$ on a standard probability space $(Y,\nu)$.
\vskip 0.02in

For every $n$-tuple $p=(p_1,p_2,..,p_n)$ as above, let $\alpha_p$ denote the diagonal product action of $\Gamma$ on $(X_p\times Y,\mu_p\times\nu)$ given by $\alpha_p(\gamma)=\beta_p(\pi(\gamma))\times \sigma(\gamma)$, for all $\gamma\in\Gamma$. 
\vskip 0.02in
 If $\alpha_p$ is orbit equivalent to $\alpha_q$,  then $h(\beta_p)=h(\beta_q)$. 
In particular, 
$\{\alpha_{(t,1-t)}\}_{t\in (0,\frac{1}{2}]}$ gives a 1-parameter family of free ergodic non-OE actions of $\Gamma$.
\endproclaim
{\it Proof.} Assume that $\alpha_p$ and $\alpha_q$ are orbit equivalent and suppose by contradiction that $h(\beta_p)<h(\beta_q)$ (after interchanging $p$ and $q$, if necessary).
By applying Theorem 4.1  we get that a quotient $\beta_p^0:\Delta\curvearrowright (X_p^0,\mu_p^0)$ of  $\beta_p$ is conjugate to the restriction of $\beta_q$ to a subgroup $\Delta_0$ of $\Delta$.  
Thus,  $$h({\beta}_p)\geq h(\beta_p^0)= h({\beta_q}_{|\Delta_0})=|\Delta/\Delta_0|h(\beta_q)\geq h(\beta_q),$$
a contradiction.
For the second assertion, note that the function $t\rightarrow h(\beta_{(t,1-t)})$ is injective on $(0,\frac{1}{2}]$.
 \hfill$\square$
\vskip 0.1in
\noindent
Since the action $\sigma:\Bbb F_n\curvearrowright (\Bbb T^2,\lambda^2)$ is free, weakly mixing, strongly ergodic and rigid, for any $2\leq n\leq \infty$ and any embedding of $\Bbb F_n$ into SL$_2
(\Bbb Z)$, we see that Theorem 5.1 implies our main result.
\vskip 0.1in
\noindent
{\bf 5.2 Final remarks.}
$(1)$. In the context from  5.1, let $M_{p}=L^{\infty}(X_p\times Y)\rtimes_{\alpha_p}\Gamma$ and assume that $\Gamma$  has Haagerup's property (e.g. $\Gamma=\Bbb F_n,2\leq n\leq \infty$).  
Then $M_p$ is not isomorphic to $M_q$, whenever $h(\beta_p)\not= h(\beta_q)$.
Indeed, following [Po06], $M_p$ is a II$_1$ factor in Popa' $\Cal H\Cal T$ class with $L^{\infty}(X_p\times Y)$ being its unique (up to conjugacy with a unitary element) HT Cartan subalgebra. Thus, isomorphism of the factors $M_p$ and $M_q$ implies orbit equivalence of the actions $\alpha_p$ and $\alpha_q$, and the claim follows from Theorem 5.1. In particular, the actions $\{\alpha_{(t,1-t)}\}_{t\in (0,\frac{1}{2}]}$ are non-von Neumann equivalent, i.e. their associated group measure space factors are non-isomorphic.
  
$(2)$. If two groups $\Gamma_1$ and $\Gamma_2$ satisfy the hypothesis of  5.1, then their product $\Gamma=\Gamma_1\times\Gamma_2$ also does. To see this, just note that if  the actions $\sigma_i:\Gamma_i\curvearrowright(Y_i,\nu_i)$   verify condition $(iii)$, then the product action $\sigma:\Gamma\curvearrowright (Y_1\times Y_2,\nu_1\times\nu_2)$   
verifies it as well. In particular, Theorem 5.1 provides uncountably many non orbit equivalent actions of $\Bbb F_m\times\Bbb F_n$, for all $2\leq m,n\leq \infty$, using an approach different from the previous ones ([MSh06], [Po08]).

$(3)$. If $\Lambda$ is an arbitrary group and $\Gamma$ is a group which satisfies the hypothesis of 5.1, then the free product $\Gamma*\Lambda$ also does. Indeed, assume that  $\sigma:\Gamma\curvearrowright (X,\mu)$ is a free,  weakly mixing, strongly ergodic, rigid action.
Following [IPP08, A.1], there exists a free action $\tilde\sigma:\Gamma*\Lambda\curvearrowright (X,\mu)$ such that $\tilde\sigma_{|\Gamma}=\sigma$.  Then it is clear that $\tilde\sigma$ is  weakly mixing and strongly ergodic. Moreover, since $\sigma$ is rigid and since $L^{\infty}(X,\mu)\subset L^{\infty}(X,\mu)\rtimes_{\sigma}\Gamma\subset L^{\infty}(X,\mu)\rtimes_{\tilde\sigma}(\Gamma*\Lambda)$, we get that $\tilde\sigma$ is a rigid action ([Po06, 4.6]). 

$(4)$. A related question (to the one considered in this paper) is to find measure preserving actions of the free groups, $\Bbb F_n$, whose associated orbit equivalence relation has trivial outer automorphism group. Note that the existence of such actions has been very recently shown by S. Popa and S. Vaes for $n=\infty$ ([PV08]) and by D. Gaboriau for $2\leq n<\infty$ ([Ga08]). 
\vskip 0.2in
\head  References\endhead

\item {[BG81]} S. I. Bezuglyi, V. Ya. Golodets: {\it Hyperfinite and II$_1$ actions for nonamenable groups}, J. Funct. Anal.  {\bf 40} (1981), no. 1, 30--44. 
\item {[Bu91]} M. Burger: {\it Kazhdan constants for SL$(3,\Bbb Z)$,}  J. Reine Angew. Math.  {\bf 413} (1991), 36--67.
\item {[CFW81]} A. Connes, J. Feldman, B. Weiss: {\it An amenable equivalence relation is generated by a single transformation},  Ergodic Theory Dynam. Systems {\bf 1} (1981), no. 4, 431--450.
\item {[CW80]} A. Connes, B. Weiss: {\it Property T and asymptotically invariant sequences}
Israel J. Math. {\bf 37} (1980), no. 3, 209–-210.
\item {[Dy59]} H. Dye: {\it On groups of measure preserving transformations I},  Amer. J. Math., {\bf 81} (1959), 119--159.
\item {[Ep07]} I. Epstein: {\it Orbit inequivalent actions of non-amenable groups}, preprint 2007,  arXiv:0707.4215. 
\item {[FM77]} J. Feldman, C.C. Moore: {\it Ergodic equivalence relations, cohomology, and von Neumann algebras. II},  Trans. Amer. Math. Soc.  {\bf 234}  (1977), no. 2, 325--359.
\item {[Ga00]} D. Gaboriau: {\it On orbit equivalence of measure preserving actions,} Rigidity in dynamics and geometry (Cambridge, 2000),  167--186, Springer, Berlin, 2002.
\item {[Ga02]} D. Gaboriau: {\it Invariants $\ell^2$ de relations d'\'equivalence et de groupes} 
Publ. math., Inst. Hautes \'Etudes Sci.,  95 (2002), no. 1, 93--150.
\item {[Ga08]} D. Gaboriau: {\it Relative Property (T) Actions and Trivial Outer Automorphism Groups}, preprint 2008, arxiv:0804.0358. 
\item {[GG88]} S.L. Gefter, V. Y. Golodets: {\it Fundamental groups for ergodic actions and actions with unit fundamental groups},  Publ. Res. Inst. Math. Sci.  {\bf 24} (1988), no. 6, 821--847.
\item {[GL07]} D. Gaboriau, R. Lyons: {\it A Measurable-Group-Theoretic Solution to von Neumann's Problem}, preprint 2007, arxiv: 0711.1643,  to appear in Invent. Math. 
\item {[GP05]} D. Gaboriau, S. Popa: {\it An uncountable family of non orbit equivalent actions of $\Bbb F_n$,}  J. Amer. Math. Soc.  {\bf 18}  (2005),  no. 3, 547--559.
\item {[Hj05]} G. Hjorth: {\it A converse to Dye's theorem}, Trans. Am. Math. Soc. {\bf 357}(2005), 3083--3103.
\item {[Ki07]} Y. Kida: {\it Classification of certain generalized Bernoulli actions of mapping class groups,} preprint 2007.
\item {[Io07a]} A. Ioana: {\it A relative version of Connes' $\chi(M)$ invariant and existence of orbit inequivalent actions},  Ergodic Theory Dynam. Systems  {\bf 27}(2007),  no. 4, 1199--1213.
\item {[Io07b]} A. Ioana: {\it Rigidity results for wreath product II$_1$ factors},  J. Funct. Anal.  {\bf 252} (2007),  no. 2, 763--791. 
\item {[Io07c]} A. Ioana: {\it Orbit inequivalent actions for groups containing a copy of $\Bbb F_2$}, preprint 2007,  arxiv: math/0701027.
\item {[IPP08]} A. Ioana, J. Peterson, S. Popa: {\it Amalgamated free products of {\it w}-rigid factors and calculation of their symmetry groups,}  Acta Math. {\bf 200} (2008), 85--153.
\item {[JSc87]} V.F.R. Jones, K. Schmidt: {\it Asymptotically invariant sequences and approximate finiteness}, 
Amer. J. Math. {\bf 109} (1987), no. 1, 91–-114.
\item {[MSh06]} N. Monod, Y. Shalom: {\it Orbit equivalence rigidity and bounded cohomology}, Ann. Math. (2) {\bf 164} (2006), 825--878. 
\item{[MvN36]} F. Murray, J. von Neumann: {\it On rings of operators}, Ann. Math. (2) {\bf 37}(1936), 116--229.   
\item {[OW80]} D. Ornstein, B. Weiss: {\it Ergodic theory of amenable group actions. I. The Rohlin lemma.}, Bull. Amer. Math. Soc. (N.S.) {\bf 2} (1980), no. 1, 161--164.
\item {[OW87]} D. Ornstein, B. Weiss: {\it Entropy and isomorphism theorems for actions of amenable groups.} J. Analyse Math. {\bf 48 }(1987), 1--141. 
\item {[Pe83]} K. Petersen: {\it Ergodic theory,} Cambridge studies in mathematics 2, 1983.
\item {[Po06]} S. Popa: {\it On a class of type II$_1$ factors with Betti numbers invariants}, Ann. Math. {\bf 163} (2006), 809--889.
\item {[Po06a]} S. Popa: {\it Strong rigidity of II$_1$ factors arising from malleable actions of $w$-rigid groups. I},  Invent. Math.  {\bf 165}  (2006),  no. 2, 369--408.
\item {[Po06b]} S. Popa: {\it Some computations of 1-cohomology groups and construction of non-orbit-equivalent actions,}  J. Inst. Math. Jussieu  {\bf 5}  (2006),  no. 2, 309--332.
\item {[Po07]} S. Popa: {\it Cocycle and orbit equivalence superrigidity for malleable actions of $w$-rigid groups},  Invent. Math. {\bf 170}  (2007),  no. 2, 243--295. 
\item {[Po07a]} S. Popa: {Deformation and rigidity for group actions and von Neumann algebras},  International
 Congress of Mathematicians. Vol. I,  445--477, Eur. Math. Soc., Z$\ddot{u}$rich, 2007. 
\item {[Po08]} S. Popa:  {\it On the superrigidity of malleable actions with spectral gap}, J. Amer. Math. Soc. {\bf 21} (2008), 981--1000.
\item {[PV08]} S. Popa, S. Vaes: {\it Actions of $\Bbb F_{\infty}$ whose II$_1$ factors and orbit equivalence relations have prescribed fundamental group}, Preprint 2008, arxiv:0803.3351. 
\item {[Sc81]} K. Schmidt: {\it Amenability, Kazhdan's property T, strong ergodicity and invariant means for ergodic group-actions},  Ergodic Theory Dynamical Systems {\bf 1}  (1981), no. 2, 223--236.
\item {[Sh05]} Y. Shalom: {\it Measurable group theory,}  European Congress of Mathematics,  391--423, Eur. Math. Soc., Z$\ddot u$rich, 2005. 
\item {[Z84]} R. Zimmer: {\it Ergodic theory and semisimple groups}
Monographs in Mathematics, 81. Birkhäuser Verlag, Basel, 1984. x+209 pp.
\enddocument